\documentclass{birkjour}
\usepackage{blindtext}
\usepackage{bm}
\usepackage{psfrag}
\usepackage[usenames,dvipsnames]{color}
\usepackage{xcolor}
\usepackage[all,cmtip,line]{xy}
\usepackage[normalem]{ulem}
\usepackage{tikz} 
\usepackage[utf8]{inputenc}
\usepackage{pgfplots}
\usepackage{scalerel}
\usepackage{amssymb}
\usepackage{caption}
\usepackage{amsmath}
\usepackage{mathabx}
\usepackage{subcaption}
\usepackage[shortlabels]{enumitem}
\usepackage[]{algorithm2e}
\usepackage{comment}

\usepackage{mathrsfs}  

\usepackage{mathtools}
\mathtoolsset{showonlyrefs}
\numberwithin{equation}{section}

\usepackage[T3,T1]{fontenc}
\DeclareSymbolFont{tipa}{T3}{cmr}{m}{n}
\DeclareMathAccent{\invbreve}{\mathalpha}{tipa}{16}

\newtheorem{theorem}{Theorem}[section]
\newtheorem{lemma}[theorem]{Lemma}
\newtheorem{corollary}[theorem]{Corollary}

\newtheorem{problem}[theorem]{Problem}
\newtheorem{definition}[theorem]{Definition}

\newtheorem{remark}[theorem]{Remark}
\newtheorem{assumption}[theorem]{Assumption}
\newtheorem{condition}[theorem]{Condition}

\newcommand{\sU}{{\normalfont\text{U}}}
\newcommand{\x}{{\bf{x}}}
\newcommand{\y}{{\bf{y}}}

\renewcommand{\div}{\operatorname{div}}

\newcommand{\cj}[1]{{\color{magenta}{#1}}}
\definecolor{forestgreen}{rgb}{0.13, 0.55, 0.13}

\newcommand{\jp}[1]{{\color{blue}#1}}

\newcommand{\todo}[1]{{\color{red}[#1]}}

\newcommand{\fh}[1]{{\color{forestgreen}{#1}}}

\newcommand{\be}{\begin{equation}}
\newcommand{\ee}{\end{equation}}

\newcommand{\cB}{\mathcal B}

\newcommand{\cO}{\mathcal O}

\newcommand{\cT}{\mathcal T}

\newcommand{\ba} {\bm a}

\newcommand{\bmu} {\bm{\mu}}

\renewcommand{\L}{\mathsf{L}}

\newcommand{\W}{{\mathsf W}}

\newcommand{\bg}{\bm{g}}

\newcommand{\norm}[1]{\left \lVert #1 \right \rVert}
\newcommand{\snorm}[1]{\left \lvert #1 \right \rvert}

\renewcommand{\div}{{\rm div}}

\newcommand{\IC}{{\mathbb C}}

\newcommand{\IN}{{\mathbb N}}

\newcommand{\IK}{{\mathbb K}}
\newcommand{\IR}{{\mathbb R}}
\newcommand{\IU}{{\mathbb U}}
\newcommand{\IT}{{\mathbb T}}
\newcommand{\IW}{{\mathbb W}}

\newcommand{\IY}{{\mathbb Y}}

\newcommand{\bp}{{\bm p}}

\newcommand{\OA}{\operatorname{\mathsf{A}}}

\newcommand{\fpcs}[1]{\operatorname{\mathsf{f}\pc}}

\newcommand{\cmspace}[3]{\mathcal{C}^{#1} \left( #2, #3 \right)}
\newcommand{\cmspaceh}[4]{\mathcal{C}^{#1,#2} \left( #3, #4 \right)}

\newcommand{\rgeoh}[2]{\mathcal{C}_b^{#1,#2}\left( (-1,1), \IR^2 \right)}

\renewcommand{\div}{\operatorname{div}}

\newcommand{\cP}{\mathcal{P}}

\newcommand{\pc}{}

\newcommand{\bla}{\boldsymbol \lambda}
\newcommand{\bphi}{\boldsymbol \phi}

\newcommand{\VH}{\bm{H}}

\newcommand{\bn}{\bm{n}}
\newcommand{\bu}{\bm{u}}
\newcommand{\bk}{\bm{k}}

\newcommand{\bz}{\bm{z}}

\newcommand{\bv}{\bm{v}}

\newcommand{\bx}{\bm{x}}
\newcommand{\by}{\bm{y}}

\newcommand{\iinterv}{\normalfont{(-1,1)}\times\normalfont{(-1,1)}}

\newcommand{\href}{\VH^{\mathrm{ref}}}

\title[Shape Holomorphy for Integral Operators on Open Arcs]{Shape Holomorphy of Boundary Integral Operators on Multiple Open Arcs}

\DeclareMathOperator{\spn}{span}

\author{Jos\'e Pinto}
\address{Facultad de Ingenier\'ia y Ciencias, Universidad Adolfo Ib\'a\~nez, Santiago, Chile}
\email{jose.pinto@uai.cl}

\author{Fernando Henr\'iquez}
\address{\'Ecole Polytechnique F\'ed\'erale de Lausanne, Switzerland}
\email{fernando.henriquez@epfl.ch}

\author{Carlos Jerez-Hanckes}
\address{Facultad de Ingenier\'ia y Ciencias, Universidad Adolfo Ib\'a\~nez, Santiago, Chile}
\email{carlos.jerez@uai.cl}


\date{December 12, 2022}

\begin{document}

\begin{abstract}
We establish shape holomorphy results for general weakly- and hyper-singular boundary integral operators arising from second-order partial differential equations in unbounded two-dimensional domains with multiple finite-length open arcs. After recasting the corresponding boundary value problems as boundary integral equations, we prove that their solutions depend holomorphically upon perturbations of the arcs' parametrizations.
These results are key to prove the shape (domain) holomorphy of the domain-to-solution maps for the associated boundary integral equations with applications in uncertainty quantification, inverse problems and deep learning.
\end{abstract}

\maketitle \footnote{This work was funded by FONDECYT INICIACION N°11230248}

\tableofcontents

\newpage

\section{Introduction}
The efficient approximation of maps with high-dimensional parametric inputs pose 
major challenges to traditional computational methods. Indeed, as the dimension of the
parametric input increases, the computational effort required to
construct surrogates of the original map may grow exponentially, thus leading to the \emph{curse of dimensionality}. In \cite{CCS15}, polynomial surrogates of high-dimensional
input maps are shown to converge independently of the dimension. Their results are derived from the
well-known one-dimensional approximation properties 
of analytic functions, a fact proven by the existence of \emph{holomorphic} extensions of the original maps onto tensor products of Bernstein ellipses in the 
complex plane. We recall that in complex variable, we refer to holomorphy as the existence of Fr\'echet derivatives of the maps on open complex subsets, which is equivalent to the existence of derivatives of arbitrary order--- analytic maps---in the same open subset. 
By varying the size of these ellipses on each parameter, namely the \emph{anisotropic} parameter dependence, one can prove
convergence rates that do not depend on the dimension of the parametric input, 
thereby \emph{breaking} the curse of dimensionality in the parametric dimension. Computationally, parametric holomorphy provides rigorous justification and construction bases for a variety of methods such as: Smolyak interpolation and quadrature \cite{ZS17,SS13};
high-order Quasi-Monte Carlo \cite{DKL14,DGGS16,DLC16,DGLS19,DGLS17}
integration (HoQMC);
deep neural network surrogates \cite{SZ19,HSZ20,OSZ22,HOS22}, 
together with its implications in Bayesian inverse problems
\cite{CS16_2,SS13,SS16}.

In this work, we consider a family of boundary value problems (BVPs) set on the complement of a finite collection of open disjoint arcs in two dimensions, with either Dirichlet or Neumann boundary conditions. We study the smoothness properties of the domain-to-solution map in the context of complex variable. Given the lack of Lipchitz regularity of the domain, traditional variational formulations cannot be applied, and hence, the existence of a holomorphic extension of the domain-to-solution map does not follow from volume-based formulations, for instance, as in \cite{CSZ18}.
Consequently, we recast the volume problems as boundary
integral equations (BIEs) posed on the collection of open arcs, as in \cite{stephane,stephan1984augmented,STE86,Stephan1987,sloan1991,JEREZHANCKES2011547,JHP20,Averseng2019,kress1996,kress2000}, and then we extend the analysis of \cite{Henriquez2021,henriquez2021_thesis} to the corresponding BIEs. 

More precisely, we will assume that each arc admits a representation arising from a suitable predefined collection of 
parametrizations. Our goal is to verify that the solutions of the BIEs depend holomorphically upon perturbations of the boundary shape. In so doing, we prove that the BIOs depend themselves holomorphically on the arcs' shape. By recalling that the inversion of linear isomorphisms defines an analytic
map, one can straightforwardly establish shape holomorphy of the domain-to-solution map as in \cite{Henriquez2021,henriquez2021_thesis}. Therein, the authors extend the solution map to the complex plane, identifying geometries with parametrizations, and prove that the corresponding map is holomorphic for Jordan arcs. Hence, by means of available complex variable results on Banach spaces, one can state that there exists a complex 
neighborhood of the collection of arc parametrizations
for which the domain-to-solution map admits a bounded holomorphic extension.
Consequently, the 
derivatives of arbitrary order of the map not only exist, the 
corresponding Taylor series expansion converges uniformly on an open neighborhood of each of the parametric arcs.

We extend significantly the work by \cite{Henriquez2021}. On one hand, our analysis encompasses a collection of open arcs as well as more general BIOs including the possibility of vector-valued ones. 
For example, we show explicit results for scalar Helmholtz and Stokes
problems, the latter also referred to as the elastic wave equation. 
The generalization is achieved by assuming the existence of a fundamental
solution with a common structure for second-order partial differential operators on two dimensions. Thus, we consider general BIOs whose kernels are given by the assumed structure of the fundamental solution, and obtain their holomorphic extensions by a
slightly more abstract version of the result for BIOs 
presented in \cite[Theorem 3.12]{Henriquez2021}. 
Also, by assuming a Maue-type representation formula, readily available for many specific operators,
one may extend the shape holomorphy of the weakly singular BIO to the hypersingular one.

In contrast with \cite{Henriquez2021}, we do not only establish the holomorphism result from  parametrizations sets to the solution in a fixed energy space, but instead consider the range of the map in a scale of functional spaces defined on the open arcs.  Indeed, in the spirit of \cite{DGLS17}, the presently obtained results
allow us to mathematically justify the use of multi-level variants of HoQMC
in the context of forward and inverse shape UQ. This provides a functional framework for high-order numerical methods, such as the one presented in \cite{JHP20}. In particular, the families of functional spaces considered here correspond to Sobolev-type spaces defined through Fourier series expansions, mapped back to open arcs by a cosine transformation (\emph{cf.}~\cite[Chapter 11]{saranen2013periodic}). 

As expected, the solutions of the Dirichlet and Neumann problems
here considered belong to highly regular spaces provided that
both the geometry and right-hand sides are regular enough.
In particular, we assume that boundary conditions are given
by the restriction of entire functions, thus having arbitrary regularity. 
Also, the arcs considered here have a limited 
regularity in a H\"older-continuous spaces. 
We thoroughly analyze how this limited smoothness
restricts the functional spaces wherein solutions of the
respective problems are sought. 

\subsection*{Outline}
The remainder of the article is as follows. After setting notation, the precise problem under consideration is given in Section \ref{sec:Problem}. Sections \ref{sec:FunctionalSpaces} through \ref{sec:hlrmfsextIKernel} introduce the necessary tools for applying abstract theorems to our particular problem.  In Section \ref{sec:abstracthlmrf} we present relevant abstract results to prove the holomorphic extension of BIOs, and also a result on how we obtain parametric holomorphism from the general holomorphic extension. In Section \ref{sec:ophlrmf} we show holomorphic extensions of the BIOs and derive the holomorphic extension of domain-to-solution maps. Similarly, in Section \ref{sec:parametrichlmrf} we conclude the parametric holomorphism of domain-to-solution maps. To illustrate our findings, in Section \ref{sec:examples} we apply these results to Helmholtz and time-harmonic elastic wave scattering by showing that the structural assumptions on the corresponding BIOs are fulfilled. Lastly, Section \ref{sec:conclu} presents conclusions and possible extensions. 

Throughout conditions will be specified. For instance, Condition \ref{condition:deltaOps} ensures that we can construct holomorphic extensions of some standard functions. In parallel, Condition \ref{condition:smh} describes the interplay among regularity parameters $s$ and $m,h$, describing functional spaces and the underlying geometry, respectively, such that the weakly- and hyper-singular operators are bounded on suitable functional spaces. These first results are the building block upon which all other conclusions are derived.
s


%

\section{Preliminaries}
Set $\imath=\sqrt{-1}$. We define the set of natural numbers $\IN$ including zero as $\IN_0 := \{ 0,1,2,\hdots\}$. Vectors and matrices are indicated by boldface symbols, while for general quantities that could be either vectors or scalars we do not use bold fonts. For a pair of vectors $\bv^1,\bv^2 \in \IC^n$, with $n \in \IN$, we define the bilinear form $ \bv^1 \cdot \bv^2 = \sum_{j=1}^n v^1_j v^2_j$, and the Euclidean norm $\|\bv^1\|^2 = \bv^1 \cdot \overline{\bv}^1$, where the conjugate of a vector is understood as component-wise conjugation. 
Also, given real numbers $a,b$, we say that $a \lesssim b $ if there exists a positive constant $c$, independent of the variables relevant for the corresponding analysis, such that $a \leq c b$. If $a \lesssim b$ and $b \lesssim a$ we write $a \cong b$.

Let $B_1,B_2$ be two Banach spaces over the field $\IK \in \{\IR,\IC\}$. 
We denote by $\mathcal{L}(B_1,B_2)$ the space of bounded linear operators
from $B_1$ to $B_2$. As it is customary, we equip it with the standard operator norm,
thus rendering it a Banach space itself.

\subsection{H\"older spaces}
\label{sec:FunctionalSpaces}

Let $\Omega_1, \Omega_2 \subset \IR^d$, $d=1,2$, be non-empty, open connected sets.
Given $m\in \IN_0$ and $\alpha \in[0,1]$, we consider the space 
$\cmspaceh{m}{\alpha}{\Omega_1}{\Omega_2}$ of functions $f:\Omega_1 \rightarrow \Omega_2$ with derivatives
up to order $m$ in $\Omega_1$ having a 
continuous extension to $\overline{\Omega}_1$, 
and such that the derivatives of order $m$ are $\alpha$-H{\"o}lder continuous.
Endowed with the norm 
\begin{align}
	\norm{f}_{\cmspaceh{m}{\alpha}{\Omega_1}{\Omega_2}}
	\coloneqq
	\sum_{\bk: |\bk| \leq m } 
	\sup_{\x \in \Omega_1}  
	\norm{\partial^{\bk} f(\bx)}
	+ 
	\sum_{\bk: \snorm{\bk} =m} 
	\sup_{\substack{\x, \y \in \Omega_1 \\ \x \neq \y}}  
	\frac{\norm{\partial^{\bk}f(\x)-\partial^{\bk}f(\y)}}{\norm{\x - \y}^\alpha},
\end{align}
where we use the standard multi-index notation \cite[page 61]{mclean2000strongly},
$\cmspaceh{m}{\alpha}{\Omega_1}{\Omega_2}$ becomes a Banach space. The case $\alpha=0$, $\mathcal{C}^{m}(\Omega_1,\Omega_2)$,  corresponds to functions
with $m$ continuous derivatives in $\overline{\Omega}_1$ with norm
\begin{align}
	\norm{f}_{\cmspaceh{m}{0}{\Omega_1}{\Omega_2}}
	\coloneqq
	\sum_{\bk: |\bk| \leq m }
	\sup_{\bx \in \Omega_1}  
	\norm{\partial^{\bk} f(\bx)}.
\end{align}
On the other hand, the case $\alpha=1$ corresponds to the one where the $m$-derivatives are Lipschitz continuous, and thus the derivatives of order
$m+1$ exist and are bounded almost everywhere (see \cite[pg. 280]{evans1998partial}).
Notice that for $m_1, m_2 \in \IN_0$ and $\alpha_1, \alpha_2 \in [0,1]$, such that $m_1 + \alpha_1 < m_2 + \alpha_2$,
one has that the inclusion $\cmspaceh{m_2}{\alpha_2}{\Omega_1}{\Omega_2} \subset \cmspaceh{m_1}{\alpha_1}{\Omega_1}{\Omega_2}$.

\subsection{Chebyshev polynomials and periodic Sobolev spaces}
Next, we recall definitions and properties of Chebyshev polynomials that will be employed to define functional spaces. For $|t|\leq1$, set $w(t) \coloneqq \sqrt{1-t^2}$, and denote by $T_n(t)$ the
$n${th} first kind Chebyshev polynomial normalized
according to 
\begin{equation}
	\int_{-1}^1 T_n(t) T_m(t) w^{-1}(t) \text{d}t 
	= 
	\delta_{n,m},
	\quad \forall \ 
	n,m \in \IN_0,
\end{equation}
being $\delta_{n,m}$ the Kronecker delta. 
Additionally, let $U_n$ denote the $n$th
Chebyshev polynomial of the second kind
normalized as follows
\begin{equation}
	\int_{-1}^1 U_n(t) U_m(t) w(t) 
	\text{d}t 
	= 
	\delta_{n,m},
	\quad \forall \ 
	n,m \in \IN_0.
\end{equation}
Furthermore, we define $\widehat{e}_n(\theta) \coloneqq \dfrac{\exp(\imath n t)}{\sqrt{2\pi}}$ as the $n$th Fourier polynomial normalized in the
$L^2(-\pi,\pi)$-norm.
For any smooth, periodic function $u :[-\pi,\pi] \rightarrow \IC$,
its Fourier coefficients are defined as
\begin{equation}
	\widetilde{u}_n 
	=
	\int_{-\pi}^\pi 
	u(\theta) 
	\widehat{e}_{-n}(\theta) 
	\text{d}\theta
	\fh{.}
\end{equation}
Similarly,  given $u: [-1,1] \rightarrow \IC$, we define two families of first kind Chebyshev
coefficients:
\begin{align}
	u_n 
	\coloneqq
	\int_{-1}^{1} 
	u(t) T_n(t) 
	\text{d}t, 
	\quad 
	\text{and} 
	\quad  
	\widehat{u}_n 
	\coloneqq
	\int_{-1}^1 
	u(t) T_n w^{-1}(t)
	\text{d}t, 
\end{align}
and two families of second kind Chebyshev coefficients:
\begin{align*}
	\ddot{u}_n 
	\coloneqq
	\int_{-1}^1 
	u(t) U_n (t)
	\text{d}t,
	\quad 
	\text{and}
	\quad 
	\widecheck{u} _n 
	\coloneqq
	\int_{-1}^1 
	u(t) U_n w(t) 
	\text{d}t
\end{align*}
Fourier coefficients of a bi-periodic
function $u:[-\pi,\pi]\times [-\pi,\pi] \rightarrow \IC$
are defined as
\begin{align}
	\widetilde{u}_{n,l} 
	\coloneqq
	\int_{-\pi}^{\pi}
	\int_{-\pi}^\pi 
	u(\theta,\phi) \widehat{e}_{-n}(\theta)\widehat{e}_{-l}(\phi) 
	\text{d}\theta 
	\text{d}\phi,
	\label{eq:biperiod}
\end{align}
and similarly for Chebyshev coefficients of bi-variate
functions on $[-1,1]\times[-1,1]$. 
We remark that the above coefficients' definitions may be 
extended to those of distributions by duality respect to the bases 
\cite[Section 5.2]{saranen2013periodic}. 

Throughout, we will make extensive use of 
periodic Sobolev spaces over $[-\pi,\pi]$,
defined for $s \in \IR$ as
\begin{equation}
	H^s[-\pi,\pi]
	\coloneqq
	\left\lbrace 
		u : 
		\norm{u}_{H^s}^2 
		= 
		\sum_{n=-\infty}^\infty (1+n^2)^s |\widetilde{u}_n|^2 
		< \infty
	\right\rbrace.
\end{equation}
We refer to \cite[Chapter 5]{saranen2013periodic} for 
a more rigorous definition.
Following \cite{Averseng2019,Jerez-Hanckes2017}, we introduce for $s\in \IR$ the following spaces defined over $(-1,1)$:
\begin{align}
	T^s 
	\coloneqq
	\left\lbrace 
		u : 
		\norm{u}_{T^s}^2 
		= 
		\sum_{n=0}^\infty (1+n^2)^s \snorm{{u}_n}^2 
		< 
		\infty
	\right\rbrace,  \\
	W^s
	\coloneqq
	\left\lbrace 
		u : 
		\norm{u}_{W^s}^2 
		= 
		\sum_{n=0}^\infty (1+n^2)^s \snorm{\widehat{u}_n}^2 
		< 
		\infty
	\right\rbrace
	\fh{.}
\end{align} 
These two functional spaces can be defined rigorously
by recalling the definition of $H^s$ and by defining two periodic
lifting operators as follows (cf.~\cite[Section 4]{Jerez-Hanckes2017})
%
%
\begin{align}
\label{eq:liffings}
(\mathcal{N}u) (\theta) := u(\cos\theta) | \sin \theta|, \quad \text{and,} \quad
(\widehat{\mathcal{N}}u)(\theta) := u (\cos\theta),
\end{align}
which again are extended to distributionsby duality along with the equivalences 
\begin{align}
\begin{split}
\label{eq:norms1}
u \in T^s \Leftrightarrow \mathcal{N}u \in H^s[-\pi,\pi], \quad \text{with,}
\quad \|u\|_{T^s} \cong \|\mathcal{N} u\|_{H^s[-\pi, \pi]}\\
 u \in W^s \Leftrightarrow \widehat{\mathcal{N}}u \in H^s[-\pi,\pi], \quad \text{with,}
\quad \|u\|_{W^s} \cong \|\widehat{\mathcal{N}} u\|_{H^s[-\pi, \pi]}.
\end{split}
\end{align}
In addition, for $s\in \IR$, we define the second kind spaces over $(-1,1)$ as
\begin{align*}
U^s := \left\lbrace u : \| u\|_{U^s}^2 = \sum_{n=0}^\infty (1+n^2)^s |\ddot{u}_n|^2 < \infty \right\rbrace, \\
Y^s := \left\lbrace u : \| {u}\|_{M^s}^2 = \sum_{n=0}^\infty (1+n^2)^s |\widecheck{u}_n|^2 < \infty \right\rbrace,
\end{align*} 
These spaces may also be
defined from periodic Sobolev spaces via  the next odd periodic liftings
\begin{align}
\label{eq:Ktheta}
(\mathcal{Z}u)(\theta) := \mathcal{N}u(\theta) \operatorname{sign}(\sin\theta), \quad \text{and,}\quad (\widehat{\mathcal{Z}}u)(\theta) := \widehat{\mathcal{N}}u(\theta) \text{sign}(\sin\theta),
\end{align}
where the sign function is defined with the convention
$\text{sign}(0)=0$.
In addition, we have the equivalences
\begin{align*}
u \in U^s \Leftrightarrow \widehat{\mathcal{Z}}u \in H^s[-\pi,\pi], \quad \text{with,}
\quad \|u\|_{U^s} \cong \|\widehat{\mathcal{Z}} u\|_{H^s[-\pi, \pi]}\\
 u \in Y^s \Leftrightarrow {\mathcal{Z}}u \in H^s[-\pi,\pi], \quad \text{with,}
\quad \|u\|_{Y^s} \cong \|\mathcal{Z} u\|_{H^s[-\pi, \pi]}.
\end{align*}
Using the previous characterizations along with, $|(\widetilde{u \sin(\cdot)})_n| \leq |\widetilde{u}_{n+1}| $, one can readily
observe that, for all $ s \in \IR$, it holds that \begin{align}\label{eq:WsYs}
	W^s \subset Y^s.
\end{align}

The dual space of $H^s$ can be identified with $H^{-s}$
in the $L^2(-\pi,\pi)$ duality pairing.
Thus, by using the lifting operators, one can identify the
dual space of $T^s$ with $W^{-s}$ and the dual of $U^s$
with $Y^{-s}$, now with respect to the $L^2(-1,1)$ duality pairing.

Furthermore, these duality identifications and \eqref{eq:WsYs} imply that
\begin{align}
\label{eq:UsTs}
U^s \subset T^s, \quad \forall \ s \in \IR.
\end{align}

From the density of the Fourier basis in $H^s$ using the inverse of the lifting operators, it is possible to deduce that the functions $\{wU_n\}_{n \in \IN}$ are dense in $U^s$. 

A mayor role in the analysis of the hyper-singular BIO is played by the mapping properties of the derivative operator (cf.~\cite[Section 4]{Jerez-Hanckes2017}). Specifically, using the density of the Chebyshev basis one can easily see that 
\begin{align}
\label{eq:devprop}
\frac{d}{dt} : U^s \rightarrow T^{s-1}, \quad \text{and,} \quad 
\frac{d}{dt} : W^s \rightarrow Y^{s-1}.
\end{align}

Also we recall from \cite[Lemma 5.3.2]{saranen2013periodic} that the periodic spaces $H^s[-\pi,\pi]$ are compactly embedded for increasing values of $s$. Hence, the same holds for the spaces $T^s,W^s,U^s,Y^s$.

Finally, depending on whether the underlying differential operator $\cP$ (see Section \ref{sec:bvproblem}) is scalar or vectorial, we use the following notation:
\begin{align*}
\IT^s = T^s,\ \IU^s = U^s, \ \IW^s = W^s, \ \IY^s = Y^s,
\end{align*}
or
\begin{align*}
\IT^s = T^s \times T^s,\ \IU^s = U^s \times U^s, \ \IW^s = W^s \times W^s, \ \IY^s = Y^s \times Y^s,
\end{align*}
respectively. We will also consider the Cartesian product spaces 
$$
\prod_{j=1}^M \IT^s, \quad \prod_{j=1}^M \IU^s,  \quad \prod_{j=1}^M \IW^s, \quad \prod_{j=1}^M \IY^s,
$$
with the standard norms.

\subsection{Bi-periodic Sobolev and H\"older spaces}

Along with the previous spaces, the forthcoming analysis will require the use of Sobolev spaces of bi-periodic functions, and an immersion result of H{\"o}lder spaces in their Sobolev counterparts. The latter will be needed in Section \ref{sec:generalbio}.

Given $s_1, s_2 $ non-negative real numbers we recall the Sobolev norm for bi-periodic functions \cite[Chapter 6]{saranen2013periodic}:
\begin{align*}
\| g\|_{s_1,s_2}^2 =
	\sum_{n=-\infty}^{\infty}
	\sum_{\ell=-\infty}^{\infty}
	(1+n^2)^{s_1} (1+\ell^2)^{s_2} 
	|\widetilde{g}_{n,l}|^2,
\end{align*}
where $\widetilde{g}_{n,l}$, are the Fourier coefficients of the bi-periodic function $g$ defined in \eqref{eq:biperiod}. Notice that in contrast to standard Sobolev spaces we could have different levels of regularity $s_1$,$s_2$ on each varible. 

\begin{lemma}
\label{lemma:cmdecay2}
Let $g \in \cmspaceh{m}{\alpha}{[-\pi,\pi]\times[-\pi,\pi]}{\IC}$
be a bi-periodic function with $m \in \IN_0$ and $\alpha \in [0,1]$.
Provided non-negative real numbers $s_1, s_2$ satisfying
$s_1+s_2 < m+\alpha$, one has that
\begin{align}\label{eq:sobolev_norm_equiv_bi_per}
	\| g\|_{s_1,s_2}^2
	\lesssim 
	\|g\|^2_{\cmspaceh{m}{\alpha}{[-\pi,\pi]\times[-\pi,\pi]}{\IC}},
\end{align}
where the implied constant is independent of $g$.
\end{lemma}

We relegate the proof of the previous lemma to Appendix \ref{appendix:lemma1}.

\begin{remark}
A similar result to Lemma \ref{lemma:cmdecay2} can be established by noticing that the smoothness of $g$ implies a certain decay speed for the Fourier coefficients. In fact, for a pair of non-negative real values $s_1,s_2$ and $g \in \cmspace{m}{[-\pi,\pi]\times[-\pi,\pi]}{\IC}$  one recovers the bound of Lemma \ref{lemma:cmdecay2}, under the more restrictive condition $s_1+s_2+1 < m$.
This requirement can be relaxed to $s_1+s_2 < m$ when considering functions in  $\cmspace{m}{[-\pi,\pi]\times[-\pi,\pi]}{\IC}$ with derivatives of order $m+1$ integrables or of bounded variation.

We also notice that, for $s_1, s_2 \in \IN_0$, one requires the slightly less restrictive condition $s_1+s_2 \leq m+\alpha$. 
\end{remark}

\subsection{Arc parametrizations}

Let $(-1,1)$ be the canonical interval. 
Throughout, we will define an open arc as the
image of a continuously differentiable,
globally invertible function 
$\bm{r} : (-1,1) \rightarrow \IR^2$.  
We further assume that the tangent vector is nowhere null
 and fix the normal vector to have the same direction of $(r'_2,-r'_1)$. 
Rigorously speaking, this definition identifies open arcs
with their corresponding parametrization instead of the set
of points that describe the arc in $\IR^2$.
The function $\bm{r}$ is also called parametrization of the arc. 

\begin{remark}
The parametrization of the open arcs 
will be taken as elements of
$\cmspaceh{m}{\alpha}{(-1,1)}{\IR^2}$.
The shape holomorphy analysis considers
\emph{complex} Banach spaces in contrast to the real-valued ones
used in open arcs parametrizations. We overcome this issue
by considering $\cmspaceh{m}{\alpha}{(-1,1)}{\IC^2}$. 
We will only use those elements with global inverse and
non null tangent vector at every point. 
The subset of such functions is denoted by $\rgeoh{m}{\alpha}$.
In some instances, we will use the Cartesian product space, 
$\prod_{j=1}^M \rgeoh{m}{\alpha}$. 
This is a subset of the product space $\prod_{j=1}^M \cmspaceh{m}{\alpha}{(-1,1)}{\IR^2}$,  equipped with the standard norm:
$$ \| \bg\|_{\prod_{j=1}^M \cmspaceh{m}{\alpha}{(-1,1)}{\IR^2}}
 = \max_{j=1,\hdots,M} \| g_j\|_{\cmspaceh{m}{\alpha}{(-1,1)}{\IR^2}}.
$$
\end{remark}

\section{Boundary Value Problems and Boundary Integral Formulation}
\label{sec:Problem}

\subsection{Boundary Value Problems on Open Arcs}
\label{sec:bvproblem}
Let us denote by $\Gamma$ the set of $M$ disjoint open finite-length arcs
$ \{ \Gamma_1, \hdots, \Gamma_M\}$, where each arc is
characterized by a parametrization
$\bm{r}_j :[-1,1] \rightarrow \Gamma_j \subset \IR^2$.
In addition, we also refer to $\Gamma$ as 
the geometric configuration of the associated problem.

Let us consider a second-order partial differential operator of the 
following general form 
\begin{align}\label{eq:partial_diff_op_P}
	\cP 
	= 
	- \sum_{j=1}^2 \sum_{k=1}^2 \partial_{x_j}{A_{j,k}} \partial_k  
	+ 
	{A},
\end{align}
where ${A_{j,k}}$ and $A$ can be constant complex-valued scalars or
or $2\times 2$ matrices. In the latter case, 
we further assume that
$\overline{\mathbf{A}}^\top_{j,k} = \mathbf{A}^\top_{k,j}$ and 
$\overline{\mathbf{A}}^\top= \mathbf{A}$.
We define the co-normal trace operator over the boundary $\Gamma_p$ as
\begin{align}
	\cB_p u 
	\coloneqq
	\sum_{j=1}^2 (\bn_p)_j \sum_{k=1}^2 {A}_{j,k} \partial_k u \big\vert_{\Gamma_p},\quad p=1,\ldots,M,
\end{align}
for any smooth function $u$ and 
where $\bn_p$ denotes the unitary normal vector of $\Gamma_p$.
Equipped with these definitions, we consider the following BVPs:
\begin{problem}[Dirichlet and Neumann BVPs]
\label{prob:BVPs}
Seek $u$ such that 
\begin{align}\label{eq:volprob}
	\cP u = 0, \quad \text{on } \ \IR^2 \setminus \overline{\Gamma},\\
	\label{eq:radcond}
	\text{Condition at infinity}(\cP),
\end{align}
complemented with either boundary conditions:  
\begin{align}
\label{eq:dircond}
u = f^D   \quad \text{on } \Gamma_p, \quad p =1,\hdots,M, \quad  \text{ (Dirichlet)},\\
\label{eq:neumanncond}
 \cB_p u = f^N  \quad \text{on } \Gamma_p, \quad p =1,\hdots,M, \quad  \text{(Neumann)}.
\end{align}
\end {problem}

Condition \eqref{eq:radcond} specifies the behavior of
$u$ far away from $\Gamma$ and it is crucial to show uniqueness. Its particular form depends on the specific partial differential operator $\cP$. Boundary data $f^D$ and $f^N$ correspond to the
right-hand sides of the Dirichlet and Neumann boundary value problems, respectively. Throughout, we assume that these are the restriction to $\Gamma$ of \emph{entire} functions in each coordinate  in $\IR^2$. Still, analytic functions on bounded domains can be used without fundamentally changing any result.

In what follows, we assume uniqueness for both Dirichlet and Neumann BVPs, while existence results will be consequence of the boundary integral formulation presented next. 

\begin{remark}
Certain assumptions on the operator $\cP$ are worth further comments. Specifically,
\begin{itemize}
	\item[(i)]
	The coefficients of 
	$\cP$ are assumed to be constants. Though this assumption is not entirely necessary, 
	one would require the coefficients
	to be at least analytic in the spatial variable. Otherwise, structural assumptions on the fundamental
	solutions of $\cP$ will not hold as detailed in the upcoming section.
	\item[(ii)]
	We have also assumed that
	$\overline{\mathbf{A}}^\top_{j,k} = \mathbf{A}^\top_{k,j}$ and 
	$\overline{\mathbf{A}}^\top= \mathbf{A}$, and also that
	$\cP$ lacks any first order derivative term. These conditions
	ensure that the co-normal trace is self-adjoint,
	thus rendering the analysis simpler.
	Yet, our results still hold without this condition by suitably modifying the associated integral operators \cite[Chapter 7]{mclean2000strongly}.
\end{itemize}
\end{remark}

\subsection{Boundary Integral Formulation}
\label{sec:bif}
We now recall the boundary integral
formulation of the BVPs for the 
partial differential operator $\cP$ introduced in 
Section \ref{sec:bvproblem}.
To this end, we assume the existence of a fundamental
solution associated to the partial differential operator $\cP$, 
which we denote by $G(\x,\y)$. 
For further details we refer to \cite[Chapter 6]{mclean2000strongly}
and references therein.
In addition, we assume that the fundamental solution
admits a decomposition of the form
\begin{align}\label{eq:funsolgen}
	G(\x,\y )  
	=
	F_1\left(\|\x -\y \| ^2\right)\log \|\x -\y \|^2
	+
	F_2\left(\|\x -\y \| ^2\right),
\end{align}
where the functions $F_1$ and 
$F_2$ are assumed to be entire complex-valued scalars or $2 \times 2$ matrices.
Furthermore, whenever $F_1$ is a scalar we assume that $F_1(0) \neq 0$ while if matrix-valued then $\mathbf{F_1}(0)$
should admit an inverse.

\begin{remark}
One could lessen the restrictions for $F_1, F_2$ and impose that they are only analytic on an open connected set of $\IC$.
If so, our results would still hold inside the analyticity domain of $F_1$, $F_2$. 
\end{remark}

Next, we introduce the single and
double layer potentials, respectively, 
on a generic open arc $\gamma$ as
\begin{align}
	\left(
		\widehat{\text{SL}}_\gamma
		\widehat{\lambda}
	\right)(\x) 
	\coloneqq
	\int_\gamma 
		G(\x,\y) 
		\widehat{\lambda}(\y)
	\text{ds}_{\y}
	\text{,}
	\ 
	\left(
		\widehat{\text{DL}}_\gamma 
		\widehat{\mu}
	\right)(\x) 
	\coloneqq
	\int_\gamma 
	\left(
		\mathcal{B}_{\bn,\y} 
		G(\x,\y)
	\right)^\top 
	\widehat{\mu}(\y) 
	\text{ds}_{\y},
\end{align} 
where $\mathcal{B}_{\bn,\by}$
denotes the co-normal trace\footnote{The general definition of the double layer potential involves the adjoint of the co-normal trace operator but under our assumptions on $\cP$ this operator is self-adjoint. } in the $\y$ variable. The densities $\widehat{\lambda}$ and $\widehat{\mu}$ are defined on
$\gamma$, and are scalar-valued or two-dimensional vectors depending on the nature of $\cP$.
The fundamental solution definition ensures that 
both potentials are homogeneous solutions of \eqref{eq:volprob} in
$\IR^2 \setminus \overline{\gamma}$.
Furthermore, we assume that both potentials
satisfy the
radiation condition \eqref{eq:radcond}. 

Let $\bm{r}: (-1,1) \rightarrow \IR^2$ be a parametrization 
of the open arc $\gamma_{\bf}$. We introduce the transformed densities: 
\begin{align*}
\lambda (\tau) := \widehat{\lambda} \circ \bm{r}(\tau)  \| \bm{r}'(\tau) \|, \quad 
\mu (\tau) := \widehat{\mu} \circ \bm{r}(\tau),\quad \tau\in(-1,1),
\end{align*}
and the pulled-back potentials
\begin{align}
	\left(
	\text{SL}_{\bm{r}} \varrho
	\right)
	(\x)
	\coloneqq
	\int_{-1}^1 
 	G(\x,\bm{r}(\tau)) \varrho(\tau)
	\text{d}\tau,
	\\
	\left(
		\text{DL}_{\bm{r}} \varrho
	\right)(\x)
	\coloneqq
	\int_{-1}^1 
	\left(\mathcal{B}_{\bn,\y}G(\x,\bm{r}(\tau))\right)^\top 
	\varrho(\tau) \|\bm{r}'(\tau)\|
	\text{d}\tau,
\end{align}
defined for $\varrho:[-1,1]\rightarrow \IC$.
From these definitions it is direct that 
$\widehat{\text{SL}}_{\gamma} \widehat{\lambda} = \text{SL}_{\bm{r}} {\lambda}$, 
and also that
$\widehat{\text{DL}}_{\gamma} \widehat{\mu} = \text{DL}_{{\bm{r}}} \mu$.

With these elements, one can now reformulate the BVPs presented in the previous section as a set of boundary integral equations (BIEs). We do so by imposing boundary conditions on indirect integral representations built via the above boundary layer potentials.

\begin{problem}[Dirichlet and Neumann BIEs]
\label{prob:BIEs}
We seek densities $\bla = (\lambda_1,\hdots,\lambda_M)$ and $\bmu = (\mu_1, \hdots, \mu_M)$, with each $\lambda_i$ and $\mu_i$ defined over $[-1,1]$ for $i = 1,\hdots,M$, such that
\begin{align}
\sum_{j=1}^M (\text{SL}_{{\bm{r}}_j} \lambda_j )\circ \bm{r}_i = f^D\circ \bm{r}_i, \quad i = 1,\hdots,M, \text{ (Dirichlet BIE)}, \\
\sum_{j=1}^M (\mathcal{B}_{\bn,\bx}\text{DL}_{{\bm{r}}_j} \mu_j )\circ \bm{r}_i = f^N\circ \bm{r}_i, \quad i = 1,\hdots,M,
\text{ (Neumann BIE)}.
\end{align}
\end{problem}
With these, we derive the following solutions for the BVPs (Problem \ref{prob:BVPs}):
\begin{align}
u=\sum_{j=1}^M \text{SL}_{{\bm{r}}_j} \lambda_j \quad \text{ (Dirichlet)}, \quad 
u=\sum_{j=1}^M \mathcal{B}_{\bn,\bx}\text{DL}_{{\bm{r}}_j} \mu_j\quad \text{ (Neumann)}.
\end{align}

We can rewrite the BIEs in matrix form:
\begin{align}
\label{eq:bios}
	\mathbf{V}_{\bm{r}_1,\dots,\bm{r}_M} \bla_{\bm{r}_1,\dots,\bm{r}_M}
	= 
	\mathbf{f}^D_{\bm{r}_1,\dots,\bm{r}_M},
	\quad 
	\mathbf{W}_{\bm{r}_1,\dots,\bm{r}_M} \bmu_{\bm{r}_1,\dots,\bm{r}_M}
	= 
	\mathbf{f}^N_{\bm{r}_1,\dots,\bm{r}_M},
\end{align}
where
\begin{equation}
	(V_{\bm{r}_1,\hdots,\bm{r}_M})_{i,j} 
	\coloneqq
	(\text{SL}_{\bm{r}_j}  )\circ \bm{r}_i
	\quad
	\text{and}
	\quad
	(W_{\bm{r}_1,\hdots,\bm{r}_M})_{i,j} 
	\coloneqq
	(\mathcal{B}_{\bn,\bx}\text{DL}_{\bm{r}_j} )\circ \bm{r}_i
\end{equation}
are weakly- and hyper-singular BIOs, and 
\begin{equation}
	(f^D_{\bm{r}_1,\hdots,\bm{r}_M})_i 
	= 
	f^D \circ \bm{r}_i
	\quad
	\text{and}
	\quad
	(f^N_{\bm{r}_1,\hdots,\bm{r}_M})_i
	\coloneqq
	f^N \circ \bm{r}_i.
\end{equation}
The weakly-singular operators can be represented
as a Lebesgue integral as follows
\begin{align*}
	(V_{\bm{r}_1,\hdots,\bm{r}_M})_{i,j}\varrho(t)
	= 
	\int_{-1}^1G(\bm{r}_i(t),\bm{r}_j(\tau)) \varrho(\tau) 
	\text{d}s,
	\quad
	t\in (-1,1),
\end{align*}
for a function $\varrho:[-1,1] \rightarrow \IR$.
On the other hand, the hyper-singular operator 
can only be expressed as a Hadamard's finite-part integral. 
However, for every $s\in \IR$, and $\varrho \in U^s$, we will assume the existence of a Maue-type representation formula of the form 
\begin{align}
\label{eq:mauesrep}
	(W_{\bm{r}_1,\hdots,\bm{r}_M})_{i,j}\varrho 
	=  
	\frac{\text{d}}{\text{d}t}
	\int_{-1}^1G(\bm{r}_i(t),\bm{r}_j(\tau)) 
	\frac{\text{d}}{\text{d}s}\varrho(\tau)
	\text{d}s 
	+ 
	\int_{-1}^1
	\widetilde{G}(\bm{r}_i(t),\bm{r}_j(\tau))\varrho(\tau) \text{d}s, 
\end{align} 
where $\widetilde{G}$ is a function with the same structure of the fundamental solution, i.e.~as in \eqref{eq:funsolgen}
. 
Such expressions for the hyper-singular operators are well known for particular cases of $\mathcal{P}$ on closed boundaries; see for example \cite{Nedelec1982,Maue1949} or the general result for scalar operators in \cite[Chapter 3.3.4]{Sauter:2011}. Yet, to the best of our knowledge, there is no known result for the general case. Similar results to the case of open arcs are derived by zero-extensions of boundary densities on the arc onto closed curves containing the arc (cf.~\cite{Jerez-Hanckes2012}).

To conclude this section, we remark that since we construct solutions of the BVPs upon  boundary potentials acting on the resulting densitites for the BIEs above, we in fact show existence results of these BVPs by showing that of the BIEs. Moreover, as it was pointed out in \cite[Remark 3.11]{JHP20}, 
the uniqueness of the boundary integral formulations can be inferred from the uniqueness of the boundary value problems, which was assumed to hold in the previous section.

\section{Holomorphic Extensions}
\label{sec:mathtools}

We now introduce the main tools to prove
the sought shape holomorphy property of the BIOs on open arcs. 
This result is stated ahead in Theorems \ref{thrm:solutionholomrext}
and \ref{thrm:bpeholomr}. The main ingredient is the holomorpic extension of BIOs to complex-valued parametrizations. In view of this, in Section \ref{sec:holom_maps_banach}
we introduce the tools concerning holomorphic
maps in general Banach spaces. In Section \ref{section:ShapeHolormsGeneral},  we consider subsets of Banach spaces that are characterized by countable many parameters (parametric subsets), and introduce the associated notion of  parametric holomorphy in Definition \ref{def:bpe_holomorphy}. We show that the general holomorpy property, previously introduced, is inherited as parametric holomorphism when a map is restricted to a parametric subset.

Based on previous work \cite{Henriquez2021,henriquez2021_thesis}, in Section \ref{sec:abstracthlmrf} we present a general framework for establishing the holomorphic extension of some general class of integral operators on arbitrary Banach spaces. Finally, in Sections \ref{sec:generalbio}, \ref{sec:funext}, and \ref{sec:hlrmfsextIKernel} we present the tools that enable us to consider the extension of the BIOs in Section \ref{sec:bvproblem} to complex parametrization in the framework of Section \ref{sec:abstracthlmrf}.


\subsection{Holomorphic Maps in Banach Spaces}\label{sec:holom_maps_banach}
Let $B_1,B_2$ be Banach spaces over the field $\IK \in \{\IR,\IC\}$.
When the underlying field of either $B_1$ or $B_2$ is $\IR$,
we say that $B_1$ or $B_2$ is a real Banach space, otherwise
we refer to them as complex Banach spaces.
One can construct a complex Banach space by taking as starting point a real one. 
For instance, suppose that $B_1$ is a real Banach space and let $\imath$ be the imaginary
unit, namely $\imath^2 =-1$. 
We define the space $B_1^{\IC}$ as the set of elements of the form $b+id$, with $b,d \in B_1$,
and we refer to it as the \emph{complexification} of $B_1$.
The complexification of $B_1$ is also a Banach space with the norm 
$\|b+id\|_{B_1^{\IC}} \coloneqq \sup_{\theta \in [0,2\pi]}\| b \cos \theta + d \sin \theta\|_{B_1}$, 
and $\IC$ as the underlying field. 
If $B_1$ is a real Banach space, given an arbitrary subset $K \subset B_1$, and $\delta$ a positive real number, we define 
\begin{align}
\label{eq:deltasets}
	K_\delta 
	= \left\lbrace k \in B_1^{\IC} \ \text{such that} \  \exists \ b \in K \ : \ \|b-k\|_{B_1^{\IC}} < \delta \right\rbrace.
\end{align}


We now introduce the notion of holomorphy in Banach spaces. 
We refer to \cite{Muj86} for a more detailed study of
complex analysis in Banach spaces.
\begin{definition}
\label{def:holmrfextension}
Let $K\subset B_1$. Assume that there exists 
an open set $\mathcal{O}\subset B^\IC_1$ such that $K \subset \mathcal{O}$.
We say that the map $f :K \subset B_1 \rightarrow B_2$ is holomorphic in
$K$ if there exists an extension of $f$ to $\mathcal{O}$, still denoted by $f$,
such that $f : \cO \rightarrow B_2$ is Fr\'echet differentiable.
\end{definition}

The next result states that if an invertible operator
admits a bounde holomorphic extension so does its inverse. 
The proof is based on well-known results
from complex variable theory. 

\begin{theorem}[{\cite[Proposition 4.20]{Henriquez2021}}]
\label{thrm:abtractinverse}
For $K\subset B_1$, let $(\OA_k)_{k\in K}$ be a family of operators in $\mathcal{L}(B_1,B_2)$ such that: 
\begin{enumerate}
	\item[(i)]
	For every $k \in K$, $\OA_k$ has a bounded inverse, i.e.
	 $\OA_k^{-1} \in \mathcal{L}({B}_2,{B}_1)$. 
	\item[(ii)]
	There exist $\delta >0$ such that the map
	${K \ni k} \mapsto \OA_k \in \mathcal{L}(B_1,{B}_2)$
	admits a bounded holomorphic extension into $K_\delta$.
\end{enumerate}
Then, there exists $\eta$, depending of $K$ and $\delta$, such that the map 
\begin{equation}
	K
	\ni
	k
	\mapsto 
	\OA_{k}^{-1} 
	\in 
	\mathcal{L}(B_2,B_1)
\end{equation}
admits a bounded holomorphic extension into $K_\eta$.
\end{theorem} 

\subsection{Parametric Holomorphy}
\label{section:ShapeHolormsGeneral}
In concrete applications, such as the ones arising
in forward and inverse computational UQ, one is usually
interested in using a parametrically defined compact set
that in turn defines the set of admissible parametric representations. 
A particular example of this construction consists in 
considering an \emph{affine}-parametric
set of the form
\begin{align}\label{eq:Kcompact}
	K 
	=
	\left\lbrace 
		k_{\by} \in B_1: 
		k_{\by} 
		= 
		k_0 
		+ 
		\sum_{n=1}^\infty y_n k_n, 
		\ 
		\by = \{y_n\}_{n \in \IN} 
		\in 
		\sU
	\right\rbrace,
\end{align}
where $k_0 \in B_1$ is fixed, $\sU\coloneqq [-1,1]^{\mathbb{N}}$
and $\{k_n\}_{n\in \IN} \subset B_1$
is  a fixed sequence, usually referred to as perturbation basis,
as one can interpret the elements $k_{\by}$ as
perturbations of the \emph{nominal} value $k_0$
modulated by the parameter sequence $\by \in \sU$.
By assuming $\{ \|k_n\|_{B_1} \} \in \ell^1(\IN)$,
one can prove that $K$ is compact in $B_1$, as shown in 
\cite[Lemma 2.7]{CD15}. Within this framework, we consider maps of the form
$F : \sU\rightarrow B_2: \by \mapsto f(k_{\by})$
where $f:B_1 \rightarrow B_2$ denotes a holomorphic map
in the sense of Definition \ref{def:holmrfextension}.
This construction renders
$F$ a Banach-space-valued map with a high-dimensional 
input.

For a rigorous study of $F$ we make use of the so-called
$(\boldsymbol{b},p,\varepsilon)$-holomorphic maps, originally introduced in
\cite{CCS15}, a key mathematical property
to break the curse of dimensionality in the approximation of parametric
maps with high-dimensional inputs. 
Specifically, for $\varrho>1$, we consider the 
Bernstein ellipse in the complex plane
\begin{align}	\mathcal{E}_{\varrho}
	\coloneqq 
	\left\{ 		\frac{\varrho+\varrho^{-1}}{2}:  z\in\IC\;\text{with} \; 1\leq \snorm{z}\leq \varrho
	\right \} 
	\subset \IC.
\end{align}
This ellipse has foci at $z=\pm 1$ and semi-axes of length 
$a\coloneqq  \frac{1}{2}({\varrho}+{\varrho}^{-1})$ and $b \coloneqq  \frac{1}{2}({\varrho}-{\varrho}^{-1})$.
Let us consider the tensorized poly-ellipse
\begin{align}
	\mathcal{E}_{\boldsymbol{\rho}} 
	\coloneqq  
	\bigotimes_{j\geq1} 
	\mathcal{E}_{\rho_j} \subset \IC^{\mathbb{N}},
\end{align}
where $\boldsymbol\rho \coloneqq  \{\rho_j\}_{j\geq1}$
is such that $\rho_j>1$, for $j\in \mathbb{N}$.
We adopt the convention $\mathcal{E}_1\coloneqq [-1,1]$
to include the case $\rho_j=1$.

\begin{definition}[{\cite[Definition 2.1]{CCS15}}]\label{def:bpe_holomorphy}
Let $X$ be a complex Banach space equipped with the norm $\norm{\cdot}_{X}$. 
For $\varepsilon>0$ and $p\in(0,1)$, we say that the map 
$$
	\sU \ni \by \mapsto u_{\by} \in X
$$
is $(\boldsymbol{b},p,\varepsilon)$-holomorphic if and only if
\begin{itemize}
	\item[(i)] 
	The map $\sU \ni \by \mapsto u_{\by} \in X$ is uniformly bounded, i.e.~
	\begin{align}
		\sup_{\by\in \sU}
		\norm{u_{\by}}_{X}
		\leq
		C_0,
	\end{align}
	for some finite constant $C_0>0$.
	\item[(ii)]
	There exists a positive sequence $\boldsymbol{b}\coloneqq \{b_j\}_{j\geq 1} \in \ell^p(\mathbb{N})$ 
	and a constant $C_\varepsilon>0$ such that for any sequence $\boldsymbol\rho\coloneqq \{\rho_j\}_{j\geq1}$ 
	of numbers strictly larger than one that is $(\boldsymbol{b},\varepsilon)$-admissible, i.e.~satisyfing
	\begin{align}
	\label{eq:admissible_polyradius}	
		\sum_{j\geq 1}(\rho_j-1) b_j 
		\leq 
		\varepsilon,
	\end{align}
	the map $\by \mapsto u_{\by}$ admits a complex
	extension $\boldsymbol{z} \mapsto u_{\boldsymbol{z}}$ 
	that is holomorphic with respect to each
	variable $z_j$ on a set of the form 
	\begin{align}
		\mathcal{O}_{\boldsymbol\rho} 
		\coloneqq  
		\displaystyle{\bigotimes_{j\geq 1}} \, \mathcal{O}_{\rho_j},
	\end{align}
	where $\mathcal{O}_{\rho_j}\subset \IC$ is an open set containing $\mathcal{E}_{\rho_j}$.
	This extension is bounded on $\mathcal{E}_{\boldsymbol\rho}$ according to
	\begin{align}
	\label{eq:bpe_hol_bound_epsilon}
		\sup_{\boldsymbol{z}\in \mathcal{E}_{\boldsymbol{\rho}}} \norm{u_{\boldsymbol{z}}}_{X}  
		\leq 
		C_\varepsilon.
	\end{align}
\end{itemize}
\end{definition}

In the context of the multiple open arcs problem
(Problem \ref{prob:BVPs}) and its
boundary integral formulation (\emph{cf.} Section
\ref{sec:bif}), we will consider a fixed geometric
configuration parametrized by functions
$\bm{r}_1^0 ,\hdots, \bm{r}_M^0$, and perturbations
defined for each individual arc in an affine manner as 
in \eqref{eq:Kcompact}.
We show that under some assumptions in the
perturbation basis, the map from the parameter space
to the boundary is holomorphic in the sense 
of Definition \ref{def:bpe_holomorphy}.

We conclude this section by a introducing a result that allows 
us to establish parametric holomorphy
in the sense of Definition \ref{def:bpe_holomorphy}.
\begin{theorem}
\label{theorem:generalbpeholomrfsm}
Assume that the sequence  $\{k_n\}_{n \in \IN} \subset B_1$ 
in \eqref{eq:Kcompact} is such that 
$\| k_n \|_{B_1} \in \ell^p(\IN)$
for some $p\in (0,1)$.
Assume that there exists $\delta>0$ such that
the map $f : K \subset B_1 \rightarrow B_2$
admits a bounded holomorphic extension
onto $K_\delta \subset B^\IC_1$.
Then, there exists $\varepsilon>0$ such that the map 
\begin{equation}
	\sU
	\ni
	\by
	\mapsto f(k_{\by}) 
	\in B_2,
\end{equation}
is $(\boldsymbol{b},p,\varepsilon)$-holomorphic with
$\boldsymbol{b} = \{\|k_n\|_{B_1}\}_{n\in \IN}$
and the same $p\in (0,1)$, and it is continuous in
the product topology.
\end{theorem}
\begin{proof}
The map $\sU \ni \by \mapsto k_{\by} \in B_1$
is $(\boldsymbol{b},p,\varepsilon)$-holomorphic. The
proof follows the same steps as that of \cite[Lemma 5.8]{CSZ18},
and we skip it for the sake of brevity. Being 
$f : K \subset B_1 \rightarrow B_2$ holomorphic itself
in the sense of Definition \ref{def:holmrfextension},
the composition of these two maps preserves this property
with the same $\boldsymbol{b} \in \ell^p(\IN)$ and $p\in (0,1)$.
The continuity statement follows by using the exact
same technique used in the proof of \cite[Lemma 5.7]{CSZ18}.
\end{proof}

\subsection{Holomorphic Integral Operators}
\label{sec:abstracthlmrf} 
We continue by following the framework introduced in \cite{Henriquez2021} so as to establish the holomorphic dependence of certain classes of BIOs.
For the sake of completeness, we revisit the most important 
results presented therein and elaborate on their extension to 
the BIOs for two-dimensional screens or cracks. 

As in Section \ref{section:ShapeHolormsGeneral}, we consider a
real-valued Banach space $B_1$, its corresponding complexification $B_1^{\IC}$
as introduced in Section \ref{sec:holom_maps_banach},
and a compact set $K \subset B_1$.
For each $k \in K$, we consider the integral operator
$\mathsf{P}_{k}$ defined as 
\begin{align}\label{eq:integral_op_Pr}
	\left(
		\mathsf{P}_{k} u
	\right)(t) 
	\coloneqq 
	\int_{-1}^{1} 
	S(t-\tau) p_{k}(t,\tau) u(\tau) 
	\text{d}s,
\end{align}
where the function $S$ does not depend on the parameter $k \in K$.
We further assume that for each $k \in K$ the integral operator $\mathsf{P}_{k}$
introduced in \eqref{eq:integral_op_Pr} defines a bounded linear 
operator between two Banach spaces
${X}$ and ${Y}$, i.e.~for each $k \in K$
we have that $\mathsf{P}_{k}  \in \mathcal{L}(X,Y)$. 
Furthermore, we assume that the continuous
functions are dense in {$X$}.
The following result enables us to construct the holomorphic
extension of the map $K\ni k \mapsto \mathsf{P}_k \in \mathcal{L}(X,Y)$. 

\begin{theorem}[{\cite[Theorem 3.12]{Henriquez2021}}] \label{thrm:abstractholm}
Assume that
\begin{enumerate}
\item[(i)]
The function $S$ is continuous everywhere except possibly at the origin.
In addition, in a neighborhood of $t=0$ we assume that
\begin{equation}
	|S(t)| \lesssim| t|^{-\beta}
\end{equation}
for some $\beta \in [0,1)$, 
and that $ p_{k} \in \cmspace{0}{(-1,1)\times(-1,1)}{\IC}$
for each $k \in K$.
\item[(ii)]
There exists a $\delta >0$ such that the map 
$$K \ni k \mapsto p_{k} \in \cmspace{0}{(-1,1)\times(-1,1)}{\IC}$$
admits a bounded holomorphic extension onto $K_\delta$ 
denoted by
\begin{equation}
	K_\delta \ni k
	\mapsto 
	p_{k,\IC} \in \cmspace{0}{(-1,1)\times(-1,1)}{\IC}.
\end{equation}
\item[(iii)] 
For $\delta$ as in item (ii), the corresponding extension of the integral operator
$\mathsf{P}_{k}$ to $K_\delta$ defined as 
\begin{align}
	(\mathsf{P}_{k,\IC}u)(t) 
	\coloneqq
	\int_{-1}^{1} S(t-\tau) p_{k,\IC}(t,\tau) u(\tau) \normalfont\text{d}s
\end{align} 
is uniformly bounded upon $K_\delta$,
i.e.~there exists a positive constant $C(K,\delta)$,
depending on $K$ and $\delta$ only, such that
\begin{align}
	\sup_{k \in K_\delta} 
	\| \mathsf{P}_{k,\IC} \|_{\mathcal{L}({X,Y})}
	<
	C(K,\delta).
\end{align}
\end{enumerate}
Then, the map 
\begin{equation}
	K
	\ni
	k
	\mapsto 
	\mathsf{P}_{k} 
	\in 
	\mathcal{L}({X,Y})
\end{equation}
admits a bounded holomorphic extension into $K_\delta$.
\end{theorem}

\begin{remark}
In \cite{Henriquez2021}, the above result is explicitly proved for continuous functions and then extended by density to an appropriate scale of Sobolev spaces. 
A close inspection of the proof reveals that this hypothesis could be further relaxed.
Indeed, it is enough to consider set of functions in $L^1(-1,1)$ that is dense in ${B_1}$ so as to apply Fubini's theorem.
\end{remark}

\subsection{BIOs on the Canonical Arc}
\label{sec:generalbio}
Based on the functional spaces defined on Section
\ref{sec:FunctionalSpaces}, we proceed to study
the mapping properties of the following three types of BIOs:
\begin{align}
	({R}_f u)(t) &\coloneqq \int_{-1}^1f(t,\tau) u(\tau) \text{d}s,\\
	({L}_fu)(t) & \coloneqq \int_{-1}^1 \log|t-\tau| f(t,\tau) u(\tau) \text{d}s,\\
	({S}_f u)(t) & \coloneqq \int_{-1}^1 \log|t-\tau| (t-\tau)^2 f(t,\tau)u(\tau)\text{d}s,
\end{align}
where $f \in \cmspaceh{m}{\alpha}{\iinterv}{\IC}$ , for some $m \in {\IN_0}$ and $\alpha \in {[0,1]}$. 
We will use the results of this section to establish the mapping properties of the weakly and hyper-singular operators, 
as this will enable us to invoke Theorem \ref{thrm:abstractholm} for the mentioned operators. The analysis follows closely that of \cite[Chapters 6 and 11]{saranen2013periodic}. In particular, we consider periodizations of the three types of integral operators, and then apply \cite[Theorem 6.1.3]{saranen2013periodic} to obtain the mapping properties of the operators. However, and in contrast to \cite[Chapter 6]{saranen2013periodic}, we consider kernels of limited regularity.

\begin{remark}
\label{rem:sigmaphi}
Let $f \in \cmspaceh{m}{\alpha}{\iinterv}{\IC}$ be given. Set
\begin{equation}
	\sigma(\theta,\phi ) 
	:= 
	f(\cos \theta, \cos \phi), \quad \text{and} 
	\quad \varphi(\theta,\phi) 
	:= 
	f(\cos \theta, \cos \phi)\sin \theta \sin \phi.
\end{equation}
One can readily observe that both functions
are bi-periodic and belong to
$\cmspaceh{m}{\alpha}{[-\pi,\pi]\times[-\pi,\pi]}{\IC}$. 
Moreover, since trigonometric functions and their derivatives
are trivially bounded, we have that
$$
\| \sigma \|_{\cmspaceh{m}{\alpha}{[-\pi,\pi]\times[-\pi,\pi]}{\IC}} \cong \| \varphi \|_{\cmspaceh{m}{\alpha}{[-\pi,\pi]\times[-\pi,\pi]}{\IC}} \cong \|f\|_{\cmspaceh{m}{\alpha}{\iinterv}{\IC}},
$$ 
where the implicit constants are independent of $f$. 
\end{remark}

\subsubsection{Operator $R_f$}
We recall the periodic lifting operators
$\widehat{\mathcal{N}}, \mathcal{N}, \mathcal{Z}, \widehat{\mathcal{Z}}$ defined
as in Section \ref{sec:FunctionalSpaces}.
When applied to the operator $R_f$ we obtain the following periodic operators: 
\begin{align*}
(\widehat{\mathcal{N}}R_f u)(\theta) &= \frac{1}{2}\int_{-\pi}^{\pi} f(\cos \theta, \cos \phi) \mathcal{N}u(\phi) \text{d} \phi,\\
(\mathcal{Z} R_f u)(\theta) &= \frac{1}{2}\int_{-\pi}^{\pi} f(\cos \theta, \cos \phi) \sin \theta \sin \phi  \widehat{\mathcal{Z}}u(\phi) \text{d} \phi.
\end{align*}
We will also make use of the following operators:
\begin{align*}
R_f^\mathcal{N} u(\theta)  &:= \frac{1}{2}\int_{-\pi}^{\pi} f(\cos \theta, \cos \phi) u(\phi) d \phi, \\
R_f^\mathcal{Z} u (\theta) &:= \frac{1}{2}\int_{-\pi}^{\pi} f(\cos \theta, \cos \phi) \sin \theta \sin \phi  u(\phi) \text{d} \phi,
\end{align*}
related to $\widehat{\mathcal{N}}R_f$, and $\mathcal{Z}R_f$ as follows
$$
\widehat{\mathcal{N}}R_f u = R_f^\mathcal{N} \mathcal{N}u , \quad \mathcal{Z}R_f u  = R^\mathcal{Z}_f \widehat{\mathcal{Z}}u.
$$
The following results follows directly from \cite[Theorem 6.1.1]{saranen2013periodic}, and Lemma \ref{lemma:cmdecay2}.
\begin{lemma}
\label{lemma:rfper}
Let $s \in \IR$, and $f \in \cmspaceh{m}{\alpha}{\iinterv}{\IC}$ the kernel function of $R_f$. If one of the following conditions is satisfied 
\begin{enumerate}
\item $s> -\frac{1}{2}$ and $s + \frac{3}{2} < m+\alpha$, or
\item $ s \leq -\frac{1}{2}$ and $-s+ \frac{1}{2} <m+\alpha$,
\end{enumerate} 
we have that $ R_f^\mathcal{N}  \in \mathcal{L}(H^{s},H^{s+1})$ and  $R_f^\mathcal{Z}  \in \mathcal{L}(H^{s},H^{s+1})$. Furthermore, they are compact operators in the corresponding spaces, and satisfy
\begin{align*}
\| R_f^\mathcal{N}\|_{ \mathcal{L}(H^{s},H^{s+1})} \lesssim \|\sigma\|_{\cmspaceh{m}{\alpha}{[-\pi,\pi]\times[-\pi,\pi]}{\IC}},\\
\| R_f^\mathcal{Z}\|_{ \mathcal{L}(H^{s},H^{s+1})} \lesssim \|\varphi\|_{\cmspaceh{m}{\alpha}{[-\pi,\pi]\times[-\pi,\pi]}{\IC}},
\end{align*}
where $\sigma, \varphi$ are defined as in Remark \ref{rem:sigmaphi}, and the unspecified constants do not depend on $f$.
\end{lemma}

\begin{proof}
Let us focus on the operator $R_f^\mathcal{N}$ as for $R_f^\mathcal{Z}$ the proof is equivalent when changing the kernel $\sigma$ with $\varphi$. We can directly see that the kernel has no singularity, and consequently, the integral operator is of arbitrary order, in particular is of order $-1$. Thus, from \cite[Theorem 6.1.3]{saranen2013periodic}, for $\nu$ arbitrary close to $\frac{1}{2}$, we have that
\begin{align*}
\begin{split}
&\| R_f^\mathcal{N}\|_{ \mathcal{L}(H^{s},H^{s+1})}^2 \lesssim \\
&\sum_{n=-\infty}^\infty \sum_{l=-\infty}^\infty  (1+n^2)^{a}(1+l^2)^{b}|\widetilde{\sigma}_{n,l}|^2+ (1+n^2)^{c}(1+l^2)^{d}|\widetilde{\sigma}_{n,l}|^2,
\end{split}
\end{align*}
where $a=s+1$, $b = \max\{\nu, |\nu-1|\}$, $c = \nu $, $d= \max \{ |s|,\nu\}$ for $s > -\frac{1}{2}$, and $a=s+1$, $b = \max\{\nu, |\nu-1|\}$, $c = \nu $, $d= \max \{ |s|,\nu\}$ for $s <-\frac{1}{2}$.

For $s=\frac{1}{2}$, we use \cite[Theorem 6.1.1]{saranen2013periodic} so as to get the same bound as before with $a =\max\{ |s+1|,\nu\}$, $b=\max\{|s|,\nu\}$, with $\nu$ as before, and parameters $c$, $d$ not specified as the second term of the right-hand side of the above inequality does not appear in this case. Hence, from Lemma \ref{lemma:cmdecay2}, if one of the conditions specified in the hypothesis holds, we obtain the bound 
$$\| R_f^\mathcal{N}\|_{ \mathcal{L}(H^{s},H^{s+1})} \lesssim \|\sigma\|_{\cmspaceh{m}{\alpha}{[-\pi,\pi]\times[-\pi,\pi]}{\IC}} < \infty,$$
and therefore, $R_f^\mathcal{N} \in \mathcal{L}(H^{s},H^{s+1})$.

Compactness follows from similar arguments as we can consider that the operator  is of order $-1-\epsilon$, for arbitrary small $\epsilon>0$, and show that $R_f^\mathcal{N} \in \mathcal{L}\left(H^s,H^{s+1+\epsilon}\right)$. Then, by the compact embedding of $H^{s}$ spaces we obtain the stated result.
\end{proof}

Using the properties of the lifting operators we derive those of the operator $R_f$. 

\begin{corollary}
\label{corollary:rfop}
Let $s  \in \IR$ be such that the hypotheses of Lemma \ref{lemma:rfper} are fulfilled. Then, we have that 
$R_f  \in \mathcal{L}(T^{s},W^{s+1})$, and $R_f \in \mathcal{L}(U^{s},Y^{s+1})$. Furthermore, they are compact operators satisfying
\begin{align*}
\| R_f\|_{ \mathcal{L}(T^{s},W^{s+1})} \lesssim \|f\|_{\cmspaceh{m}{\alpha}{\iinterv}{\IC}},\\
\| R_f\|_{ \mathcal{L}(U^{s},Y^{s+1})} \lesssim \|f\|_{\cmspaceh{m}{\alpha}{\iinterv}{\IC}},
\end{align*}
with implicit constants independent of $f$.
\end{corollary}
\begin{proof}
From the properties of the periodic lifting operator \eqref{eq:norms1}, and the previous Lemma, one deduces that
\begin{align*}
\| R_f u\|_{W^{s+1}} \cong \| \widehat{\mathcal{N}} R_f u \|_{H^{s+1}} = \|R_f^\mathcal{N} \mathcal{N}u \|_{H^{s+1}} \leq \| R_f^\mathcal{N}\|_{\mathcal{L}(H^{s},H^{s+1})} \|\mathcal{N} u \|_{H^{s}}. 
\end{align*}
We recall that $\|\mathcal{N} u \|_{H^{s}} \cong \|u\|_{T^{s}}$ (cf.~Equation \eqref{eq:norms1}), and also by the previous Lemma and Remark \ref{rem:sigmaphi}, we get 
\begin{align*}
\| R_f^\mathcal{N}\|_{\mathcal{L}(H^{s},H^{s+1})} \lesssim \|f\|_{\cmspaceh{m}{\alpha}{\iinterv}{\IC}},
\end{align*}
as stated. The proof is analogous for $U^{s},Y^{s+1}$, but using $R^{\mathcal{Z}}_f$ instead of $R^{\mathcal{N}}_f$.
\end{proof}
\subsubsection{Operator $L_f$}
As in the previous case, we consider the two lifting versions of $L_f$:
\begin{align}
\label{eq:Lsplit}
\begin{split}
(\widehat{\mathcal{N}}L_fu)(\theta) &= \frac{\log{2}}{2} \int_{-\pi}^{\pi} f(\cos \theta , \cos \phi) \mathcal{N}u(\phi) \text{d}\phi  \\&+ \int_{-\pi}^{\pi} f(\cos \theta, \cos \phi) \log \left\vert \sin \left(\frac{\theta-\phi}{2}\right) \right\vert \mathcal{N}u(\phi) \text{d}\phi,
\end{split}
\end{align}
and
\begin{align}
\label{eq:Lsplitodd}
\begin{split}
(\mathcal{Z}L_fu)(\theta) &= \frac{\log{2}}{2} \int_{-\pi}^{\pi} f(\cos \theta , \cos \phi)\sin \theta \sin \phi \widehat{\mathcal{Z}}u(\phi) \text{d}\phi \\&+ \int_{-\pi}^{\pi} f(\cos \theta, \cos \phi)\sin \theta \sin \phi \log \left\vert \sin \left(\frac{\theta-\phi}{2}\right) \right\vert \widehat{\mathcal{Z}}u(\phi) \text{d}\phi.
\end{split}
\end{align}
We see that these two operators can be characterized as the sum of a regular operator plus a logarithmic one. The logarithmic part gives rise to an operator of order $-1$. Thus, by the same arguments used in the analysis of $R_f$, we arrive at the following result:

\begin{corollary}
\label{corollary:lfop}
For $s \in \IR$, let the hypotheses of Lemma \ref{lemma:rfper} hold. Then,
 $ L_f  \in \mathcal{L}(T^{s},W^{s+1})$ and  $L_f  \in \mathcal{L}(U^{s},Y^{s+1})$. Furthermore, the bounds
\begin{align*}
\| L_f\|_{ \mathcal{L}(T^{s},W^{s+1})} \lesssim \|f\|_{\cmspaceh{m}{\alpha}{\iinterv}{\IC}},\\
\| L_f\|_{ \mathcal{L}(U^{s},Y^{s+1})} \lesssim \|f\|_{\cmspaceh{m}{\alpha}{\iinterv}{\IC}},
\end{align*}
hold with unspecified constants independent of $f$.
\end{corollary}
\subsubsection{Operator $S_f$}
Finally, consider the $S_f$ operator, whose periodic liftings are 
\begin{align*}
	&(\widehat{\mathcal{N}}S_fu)(\theta) 
	= 
	\frac{\log 2}{2}\int_{-\pi}^\pi
(\cos \theta - \cos \phi)^2 f(\cos \theta, \cos \phi)\mathcal{N}u(\phi) \text{d}\phi \\
&+4 \int_{-\pi}^\pi \log \left\vert \sin \left(\frac{\theta-\phi}{2}\right) \right\vert \sin^2\left( \frac{\theta-\phi}{2} \right)
\sin^2\left( \frac{\theta+\phi}{2} \right) f(\cos \theta, \cos \phi) \mathcal{N}u (\phi) \text{d}\phi, 
\end{align*} 
and
\begin{align*}
	&(\mathcal{Z} S_fu)(\theta) 
	= 
	\frac{\log 2}{2}
	\int_{-\pi}^\pi
	(\cos \theta - \cos \phi)^2 
	f(\cos \theta, \cos \phi)
	\sin \theta \sin \phi \widehat{\mathcal{Z}}u(\phi) 
	\text{d}\phi+\\
	&	
	4 \int_{-\pi}^\pi 
	\log \left\vert \sin \left(\frac{\theta-\phi}{2} \right)\right\vert 
	\sin^2\left( \frac{\theta-\phi}{2} \right)
	\sin^2\left( \frac{\theta+\phi}{2} \right) 
	f(\cos \theta, \cos \phi) 
	\sin \theta 
	\sin \phi \widehat{\mathcal{Z}}u (\phi) \text{d} \phi.
\end{align*}
While these operators are of order $-3$, we will consider them as a compact operator of order $-1$. This can be done by analyzing the mapping properties from $T^s$ (resp.~$U^s$) to $W^{s+1+\epsilon}$ (resp.~$Y^{s+1+\epsilon}$). In particular, we can select $\epsilon$ small enough such that the same conditions of Lemma \ref{lemma:rfper} apply, hence we obtain the following result. 

\begin{corollary}
\label{corollary:sfop}
Let $s \in \IR$, such that the same conditions of Corollary \ref{corollary:lfop} are satisfied. Then,
 $ S_f  \in \mathcal{L}(T^{s},W^{s+1})$ and  $S_f  \in \mathcal{L}(U^{s},Y^{s+1})$, and it is compact in both cases. Moreover, we have the bounds: 
\begin{align*}
\| S_f\|_{ \mathcal{L}(T^{s},W^{s+1})} \lesssim \|f\|_{\cmspaceh{m}{\alpha}{\iinterv}{\IC}},\\
\| S_f\|_{ \mathcal{L}(U^{s},Y^{s+1})} \lesssim \|f\|_{\cmspaceh{m}{\alpha}{\iinterv}{\IC}},
\end{align*}
with unspecified constants independent of $f$.
\end{corollary}
\subsection{Holomorphic Functions}
\label{sec:funext}
In the ensuing analysis, we show the existence of 
holomorphic extensions for certain 
recurrently functions (cf.~\cite[Section 4.1]{Henriquez2021}). 
However, as we are working with open arcs, the functions
considered herein are not periodic. 
The analysis provided in this section
lies in the context of spaces
of the form $\mathcal{C}^{m,\alpha}((-1,1),\IR^2)$,
with $m \in \IN$ and $\alpha \in [0,1]$, 
with at least $m+\alpha>2$,
instead of 
twice continuously differentiable, periodic functions. Consequently, the holomorphic extension of functions for multiples
arcs requires suitable sets such as the ones below.

\begin{definition}\label{def:single_admissible_arc_param}
Let $m \in \IN$ and $\alpha \in [0,1]$.
We say that $K$ is an $(m,\alpha)$-admissible
set of arc parametrizations if $K\subset \rgeoh{m}{\alpha}$
and if $K$ is a compact subset of
$\mathcal{C}^{m,\alpha}((-1,1),\IR^2)$.
\end{definition}

When dealing with multiple arcs we further need to impose
that no two pair of arcs intersect or touch each other. 
In the following, we work under the assumption stated below.

\begin{assumption}\label{assump:admissible_arc_param}
Let $K^1, \hdots,K^M$ be a collection of $M \in \IN$ 
$(m,\alpha)$-admissible set of arc parametrizations,
in the sense of Defintion \ref{def:single_admissible_arc_param},
for some $m \in \IN$ and $\alpha \in [0,1]$.
For each $i,j \in \{1,\hdots,n\}$ with $i\neq j$ it holds
$$
	\inf_{(\bm{r},\bp) \in K^i \times K^j}
	\inf_{(t,\tau) \in (-1,1)\times(-1,1)}
	\norm{\bm{r}(t) - \bp(\tau)}
	>0.
$$
\end{assumption}

Due to the structure of the $G(\cdot,\cdot)$, previously introduced in 
\eqref{eq:funsolgen}, we extensively make use of the logarithmic function,
which admits an holomorphic extension in $\IC \setminus (-\infty,0]$. 
Similarly, we also use the holomorphic extension of the squared distance function between two points located in two---not necessarily different---arcs. 

For $\bm{r}, \bp:(-1,1) \rightarrow \IR^2$, the squared distance is defined as $d_{\bm{r},\bp}^2(t,\tau) = \| \bm{r}(t)-\bp(\tau)\|^2$, and its extension to complex parametrizations takes the form (cp.~\cite[Section 4.1]{Henriquez2021})
\begin{align}
	d^2_{\bm{r},\bp}(t,\tau) = {(\bm{r}(t)-\bp(\tau))\cdot(\bm{r}({t})-\bp(\tau))},
\end{align}
where we have used the Euclidean inner product in the bilinear sense,
as the inclusion of the complex conjugation prevents the existence of any
holomorphic extension. 
Whenever $\bm{r} = \bp$ we use the notation $d^2_{\bm{r}} = d^2_{\bm{r},\bm{r}}$. 

The following conditions will be later required to establish bounds on how large are the regions where the BIOs admit holomorphic extensions.

\begin{condition}
\label{condition:Qself}
Let $m\in \IN$ and $\alpha\in [0,1]$, and
let $K$ be an $(m,\alpha)$-admissible set of arc parametrizations.
The value $\delta>0$ satisfies
$$ 
	\delta
	<
	\sqrt{\mathcal{I}_Q^2+\mathcal{S}^2_Q}-\mathcal{S}_Q,
$$
where
\begin{align*}
\mathcal{I}_Q:= \inf_{\bm{r} \in K} \inf_{t \in (-1,1)} \| \bm{r}'(t)\|,
\quad
\text{and}
\quad
\mathcal{S}_Q:= \sup_{\bm{r} \in K} \sup_{t \in (-1,1)} \|\bm{r}'(t)\|.
\end{align*} 
\end{condition}

\begin{condition}
\label{condition:dcross}
Let $m\in \IN$ and $\alpha\in [0,1]$, 
and let $K^1,K^2$ be two $(m,\alpha)$-admissible sets of arc parametrizations satisfying
Assumption \ref{assump:admissible_arc_param}.
The values $\delta_1,\delta_2>0$ satisfy
$$
\delta_1 + \delta_2 < \sqrt{\mathcal{I}_d
^2+\mathcal{S}^2_d}-\mathcal{S}_d,
$$
where 
\begin{align*}
\mathcal{I}_d&:=\inf_{(\bm{r},\bp) \in K^1 \times K^2} \inf_{(t,\tau) \in (-1,1)\times(-1,1)}
 \| \bm{r}(t) - \bp(\tau) \|,\\
\mathcal{S}_d&:= \sup_{\bm{r} \in K^1_{\delta_1}} \sup_{t \in (-1,1)} \| \bm{r}(t)\| +
\sup_{\bp \in K^2_{\delta_2}} \sup_{t \in (-1,1)} \| \bp(t)\|. 
\end{align*}
\end{condition}


The following result ensures that the square of the
distance function admits a bounded holomorphic extension
to a set of the form $K_\delta$ for some $\delta>0$.

\begin{lemma} 
\label{lemma:squaredistance}
Let $m\in \IN$ and $\alpha\in [0,1]$, 
and let $K^1,K^2$ be two $(m,\alpha)$-admissible sets of arc parametrizations.
\begin{itemize}
\item[(i)]
For any pair $\tau_1,\tau_2 > 0$, the map 
\begin{equation}
	K^1 \times K^2 
	\ni
	(\bm{r},\bp) 
	\mapsto 
	d_{\bm{r},\bp}^2 
	\in 
	\cmspaceh{m}{\alpha}{\iinterv}{\IC}
\end{equation}
admits a bounded holomorphic extension into
$K^1_{{\tau}_1} \times K^2_{{\tau}_2}$.
\item[(ii)]
There exist $\delta_1>0$ and $\delta_2>0$
satisfying Condition \ref{condition:dcross}
and $\eta>0$ such that 
\begin{align}
	\inf_{(\bm{r},\bp) \in K^1_{\delta_1}  \times K^2_{\delta_2}} 
	\inf_{(t,\tau)\in \iinterv}
	\mathfrak{Re}\{d_{\bm{r},\bp}^2(t,\tau)\}
	\geq \eta
	>
	0.
\end{align}
\end{itemize}
\end{lemma}
\begin{proof}
See Appendix \ref{proof:lem1}.
\end{proof}

Another relevant function required
to establish holomorphic extensions of our
integral operators is the following.
For each arc parametrization $\bm{r}: (-1,1)\rightarrow \IR^2$
we define $Q_{\bm{r}}:(-1,1) \times (-1,1) \rightarrow \IR$ as
\begin{align}
	Q_{\bm{r}}(t,\tau) 
	\coloneqq
	\begin{cases}
		\dfrac{d^2_{\bm{r}}(t,\tau)}{(t-\tau)^2}, \quad t\neq \tau, \\
		\bm{r} '(t) \cdot \bm{r} '(t), \quad t =s,
	\end{cases}
	\quad
			(t,\tau) \in (-1,1) \times (-1,1).
\end{align}
One can prove the following result.
\begin{lemma}
\label{lemma:Qfun}
Let $m\in \IN$ and $\alpha\in [0,1]$, and
let $K$ be an $(m,\alpha)$-admissible set
of arc parametrizations.
\begin{itemize}
\item[(i)]
There exists $\delta= \delta(K)>0$, depending upon 
$K$ only and satisfying Condition \ref{condition:Qself}, such that the map
\begin{align*}
	K
	\ni
	\bm{r}
	\mapsto
	Q_{\bm{r}} 
	\in 
	\cmspaceh{m-1}{\alpha}{\normalfont{(-1,1)\times(-1,1)}}{\IC}
\end{align*}
admits a bounded holomorphic extension into
$K_\delta$.
\item[(ii)]There exists a constant 
$\zeta= \zeta(K,\delta)>0$, depending upon $K$
and $\delta$ only, such that
\begin{align}
	\inf_{\bm{r} \in K_\delta} 
	\inf_{(t,\tau) \in \normalfont{(-1,1)\times(-1,1)}} 
	\mathfrak{Re}\{Q_{\bm{r}}(t,\tau)\} \geq\zeta, \\
	\inf_{\bm{r} \in K_\delta} 
	\inf_{(t,\tau) \in \normalfont{(-1,1)\times(-1,1)}}
	\mathfrak{Re}\{Q^{-1}_{\bm{r}}(t,\tau)\}\ \geq\zeta,
\end{align}
where $Q^{-1}_{\bm{r}}$ represents the multiplicative inverse $1/Q_{\bm{r}}$.
\item[(iii)]
The map
\begin{align*}
	K
	\ni
	\bm{r} 
	\mapsto
	Q^{-1}_{\bm{r}} 
	\in 
	\cmspaceh{m-1}{\alpha}{\normalfont{(-1,1)\times(-1,1)}}{\IC}
\end{align*} 
admits a bounded holomorphic extension into $K_{\delta}$. 
\end{itemize}
\end{lemma} 
\begin{proof}
See Appendix \ref{proof:lem2}.
\end{proof}


\subsection{Holomorphic Extension of Integral Kernels}
\label{sec:hlrmfsextIKernel}
We now show that the kernels of the weakly- and hyper-singular BIOs---according to the representations in \eqref{eq:funsolgen}--- have holomorphic extensions. We do so by extending the integral kernels using our previous results on the functions $d_{\bm{r},\bp}^2$ and $Q_{\bm{r}}$.

By \eqref{eq:funsolgen}, for
two---not necessarily different---arc parametrizations $\bm{r},\bp: (-1,1)\rightarrow \IR^2$, one can write
\begin{equation}
	\label{eq:G_disjoint_arcs}
	G(\bm{r}(t),\bp(\tau)) 
	= 
	\log(d^2_{\bm{r},\bp}(t,\tau)) F_1(d^2_{\bm{r},\bp}(t,\tau))
	+
	F_2(d^2_{\bm{r},\bp}(t,\tau)).
\end{equation}
The next result follows straightforwardly from Lemma \ref{lemma:squaredistance}
and the structure assumed for $G(\cdot,\cdot)$ in \eqref{eq:funsolgen}.
\begin{lemma}
\label{lemma:crosskernel}
Let $m\in \IN$ and $\alpha\in [0,1]$, 
and let $K^1,K^2$ be two $(m,\alpha)$-admissible sets of arc parametrizations satisfying
Assumption \ref{assump:admissible_arc_param}.
Then, there exist $\delta_1,\delta_2>0$ 
satisfying Condition \ref{condition:dcross}
such that the map
\begin{equation}
	\label{eq:analytic_G_r_p}
	K^1 \times K^2 
	\ni
	(\bm{r},\bp) 
	\mapsto 
	G(\bm{r}(t),\bp(\tau)) 
	\in \cmspaceh{m}{\alpha}{\iinterv}{\IC}
\end{equation}
admits a bounded holomorphic extension into
$K^1_{\delta_1} \times K^2_{\delta_2}$.
\end{lemma}
\begin{proof}
By Lemma \ref{lemma:squaredistance} (ii), 
there exist $\delta_1,\delta_2>0$ such that the real part
of the logarithm argument
in \eqref{eq:G_disjoint_arcs} is bounded
from below away from zero. Hence, the function 
$G(\bm{r}(t),\bp(\tau))$ is well defined for all $(t,\tau) \in (-1,1) \times (-1,1)$,
and for any non-intersecting arc parametrizations
$\bm{r},\bp: (-1,1) \rightarrow \IR^2$. Furthermore,
by Lemma \ref{lemma:squaredistance} (i) 
along with the fact that the logarithm is analytic in
the branch cut $\IC \backslash(-\infty,0]$, one concludes that the map
in \eqref{eq:analytic_G_r_p} admits a bounded holomorphic extension in
$K^1_{\delta_1} \times K^2_{\delta_2}$.
\end{proof}

For the self-interaction case, i.e.~$\bm{r}=\bp$, the result does not
follow from the arguments used in the proof
of Lemma \ref{lemma:crosskernel}, as a logarithmic singularity
inevitably occurs at $d_{\bm{r}}(t,t)=0$, thus breaking the analyticity of
the logarithmic.
In this case, we consider the following decomposition:
\begin{align}
\label{eq:kernelsplit}
	G(\bm{r}(t),\bm{r}(\tau) ) 
	= 
	G_{\bm{r}}^R(t,\tau)+G_{\bm{r}}^S(t,\tau),
\end{align}
where 
\begin{align}
	\label{eq:GR}
	G_{\bm{r}}^R(t,\tau) &
	\coloneqq
	\left(\log Q_{\bm{r}}(t,\tau)\right) F_1(d^2_{\bm{r}}(t,\tau))+F_2(d^2_{\bm{r}}(t,\tau)),\\
	\label{eq:GS}
	G_{\bm{r}}^S(t,\tau) &
	\coloneqq
	2\log|t-\tau| F_1(d^2_{\bm{r}}(t,\tau)).
\end{align}
Notice that now the logarithmic singularity has been isolated in
the term $G_{\bm{r}}^S$ defined in \eqref{eq:GS}.
Furthermore, it does not depend 
on any arc parametrization, and we have following result.
\begin{lemma}
\label{lemma:selfkernel}
Let $m\in \IN$ and $\alpha\in [0,1]$, and
let $K$ be an $(m,\alpha)$-admissible set of
arc parametrizations.
Then there exists $\delta>0$ satisfying Condition \ref{condition:Qself}
such that
\begin{align*}
	K \ni \bm{r} &\mapsto G_{\bm{r}}^R
	\in 
	\cmspaceh{m-1}{\alpha}{\iinterv}{\IC}\text{ and}\\
	K \ni \bm{r} &\mapsto F_1(d^2_{\bm{r}})
	\in 
	\cmspaceh{m}{\alpha}{\iinterv}{\IC}
\end{align*}
admit bounded holomorphic extensions onto $K_\delta$.
\end{lemma}
\begin{proof}
The only part that does not follow directly is the logarithmic term of $G_{\bm{r}}^R$. However, by Lemma \ref{lemma:Qfun}, we are again in the holomorphic domain of the logarithmic and one retrieves the above statements. 
\end{proof}

\begin{remark}
\label{remark:lessregularitykernel}
We have assumed that the functions $F_1$ and $F_2$ 
in the decomposition of $G(\bx,\by)$ stated in \eqref{eq:funsolgen}
are entire, and that they depend solely on the square of the distance
between points $\bx$ and $\by$.
However, less restrictive cases are to be considered.
For example, we will consider cases where
$F_1$, $F_2$ 
are replaced by the functions $G_1, G_2$ that
take the following form:  
\begin{equation}
	G_j(t,\tau) 
	= 
	f_j
	\left(
		\bm{r}'(t),\bp'(\tau))F_j(d_{\bm{r},\bp}^2(t-\tau)
	\right), 
	\quad j=1,2,
\end{equation}
where $f_j$ are entire in both coordinates, and $F_j$ are as before for $j=1,2$. Under this assumption, both Lemmas \ref{lemma:crosskernel} and \ref{lemma:selfkernel} still hold true but the space $\cmspaceh{m}{\alpha}{\iinterv}{\IC}$ has to be replaced by $\cmspaceh{m-1}{\alpha}{\iinterv}{\IC}$,
as the functions $f_1,$ and $f_2$ now depend on the derivative of the arc parametrizations.
This loss of one order of regularity has no effect, as we need to consider the holomorphic extension of the full kernel function which, by the first map in Lemma \ref{lemma:selfkernel}, needs to lie in $\cmspaceh{m-1}{\alpha}{\iinterv}{\IC}$ anyway.
\end{remark}

\section{Shape Holomorphy of Domain-to-Solution Maps}
\label{sec:ophlrmf}

We now study the holomorphic properties of the boundary-to-solution maps:
\begin{align}
	(\bm{r}_1,\hdots,\bm{r}_M) 
	\mapsto 
	\bla_{\bm{r}_1,\hdots,\bm{r}_M}
	\quad
	\text{and}
	\quad
	(\bm{r}_1,\hdots,\bm{r}_M) 
	\mapsto 
	\bmu_{\bm{r}_1,\hdots,\bm{r}_M},
\end{align}
where $\bla_{\bm{r}_1,\hdots,\bm{r}_M},\bmu_{\bm{r}_1,\hdots,\bm{r}_M}$
are the solutions of the Dirichlet
and Neumann boundary integral formulations introduced
in \eqref{eq:bios}. The study is carried out in three main steps:
\begin{itemize}
    \item[(i)] By Theorem \ref{thrm:abstractholm} and the results from Sections \ref{sec:generalbio}, \ref{sec:funext}, and \ref{sec:hlrmfsextIKernel}, we show that the following maps have holomorphic extensions 
$$
(\bm{r}_1,\hdots,\bm{r}_M) 
	\mapsto \mathbf{V}_{\bm{r}_1,\hdots,\bm{r}_M}, \quad (\bm{r}_1,\hdots,\bm{r}_M) 
	\mapsto \mathbf{W}_{\bm{r}_1,\hdots,\bm{r}_M},
$$
on proper compact subsets of $\prod_{j=1}^M \cmspaceh{m}{\alpha}{(-1,1)}{\IR^2}$, where the boundary integral operatos $\mathbf{V}_{\bm{r}_1,\hdots,\bm{r}_M}$, $\mathbf{W}_{\bm{r}_1,\hdots,\bm{r}_M}$ are those from Section \ref{sec:bvproblem}. 
\item[(ii)] We show that the previous operators have inverses, and then use Theorem \ref{thrm:abtractinverse} to obtain the holomorphic extensions of the boundary to solution map. 
\item[(iii)] We consider arc parametrizations determined by a countable set of parameters and we study the parametic holomorphism of the domain-to-solution map. We do so by combining the above results and the abstract ones in Section \ref{section:ShapeHolormsGeneral}.
\end{itemize}
Stages one and two are carried out in Sections \ref{ref:single_interaction}, for a single arc, and \ref{sec:multiple_arcs} for multiples arcs. The final step is presented in Section \ref{sec:parametrichlmrf}.
\subsection{Single Interaction}
\label{ref:single_interaction}
Firstly, let us study the weakly-singular BIO
between two arc parame-trizations.
Let $\bm{r},\bp \in \rgeoh{m}{\alpha}$
be the parametrization of two open arcs. 
For $u$ defined in $[-1,1]$ we set
\begin{equation}\label{eq:weakly_operator_V}
	\left(V_{\bm{r},\bp}u\right)(t)
	= 
	\int_{-1}^1G(\bm{r}(t),\bp(\tau)) u(\tau) 
	\text{d}\tau,
	\quad
	t \in (-1,1).
\end{equation}
Following the notation of Section \ref{sec:bif},
we have $V_{\bm{r}_i, \bm{r}_j} = (V_{\bm{r}_1,\hdots,\bm{r}_M})_{i,j}$,
for $i,j=1,\dots,M$.

Due to the fundamental solution structure \eqref{eq:funsolgen} as well as Lemmas \ref{lemma:crosskernel}
and \ref{lemma:selfkernel}, the operator $V_{\bm{r},\bp}$
\eqref{eq:weakly_operator_V} can be expressed in terms of $R_f, L_f$ and $S_f$ introduced in Section \ref{sec:generalbio}.
The function $G(\bm{r}(t),\bp(\tau))$ and its suitable
decomposition will play the role of
$f$ in the aforementioned operators on the canonical arc.
This analysis in performed thoroughly in Lemma \ref{lemma:vcomp} 
ahead. For the hyper-singular BIO, we assume the existence of a suitable Maue's formula 
so as to reuse the shape holomorphy result for BIOs resembling weakly-singular BIOs.

Based on results of Section \ref{sec:generalbio},
we introduce the following condition that will enable
us to ensure the continuity of the integral operators
on appropriate spaces.



\begin{condition}
\label{condition:smh}
For $m \in \IN$, $\alpha \in [0,1]$, $s \in \IR$,
one of the following conditions needs to hold:
\begin{enumerate}
\item[(i)] 
$s > -\frac{1}{2}$, and $s+
\frac{5}{2}<m+\alpha$,
\item[(ii)] 
$ s \leq -\frac{1}{2}$, and $\frac{3}{2} -s < m+\alpha$.
\end{enumerate}
\end{condition}

Observe that these conditions are exactly as those required in
Lemma \ref{lemma:rfper} but with $m-1$ instead of $m$.
This is due to the loss of regularity in the kernel with
respect to the parametrization.
Equipped with these results, we can state the main result concerning the holomorphic dependence of the operator $V_{\bm{r},\bp}$ upon a set of arc parametrizations.
\begin{lemma}\label{lemma:vcomp} 
Assume that Condition \ref{condition:smh} holds for some
$m \in \IN$, $\alpha \in [0,1]$, and $s \in \IR$.
\begin{enumerate}
	\item[(i)] 
	Let $K$ be an $(m,\alpha)$-admissible set
	of arc parametrizations. Then, there exists $\delta>0$,
	depending only on $K$ and
	satisfying Condition \ref{condition:Qself},
	such that
	\begin{align}
		K
		\ni
		\bm{r}
		\mapsto V_{\bm{r},\bm{r}} 
		\in 
		\mathcal{L}(\mathbb{T}^s, \mathbb{W}^{s+1})
		\quad
		\text{and}
		\quad
		K
		\ni 
		\bm{r}
		\mapsto V_{\bm{r},\bm{r}} 
		\in 
		\mathcal{L}(\mathbb{U}^s, \mathbb{Y}^{s+1})
\end{align}
admit bounded holomorphic extensions into $K_\delta$.
Furthermore, for each $\bm{r} \in K_{\delta}$ 
it holds that $V_{\bm{r},\bm{r}} \in \mathcal{L}(\IT^s,\IW^{s+1})$
is a Fredholm operator of index zero.
	\item[(ii)]  
	Let $K^1,K^2$ be two $(m,\alpha)$-admissible 
	sets of arc parametrizations satisfying
	Assumption \ref{assump:admissible_arc_param}.
	Then, there exist $\delta_1$, $\delta_2>0$, 
	depending on $K^1$ and $K^2$ and satisfying Condition \ref{condition:dcross},
	such that the maps
	 \begin{align}
	 	K^1 \times K^2 
		&
		\ni
		(\bm{r},\bp)
			\mapsto V_{\bm{r},\bp} 
		\in 
		\mathcal{L}(\mathbb{T}^s,\mathbb{W}^{s+1})
		\quad
		\text{and}\\
 		K^1 \times K^2
		&
		\ni
		(\bm{r},\bp)
			\mapsto V_{\bm{r},\bp} 
		\in 
		\mathcal{L}(\mathbb{U}^s, \mathbb{Y}^{s+1})
	\end{align}
	admit bounded holomorphic extensions into
	$K^1_{\delta_1} \times K^2_{\delta_2}$.
	Moreover, for any 
	$(\bm{r},\bp) \in K^1_{\delta_1} \times K^2_{\delta_2}$
	the maps $V_{\bm{r},\bp} \in \mathcal{L}(\mathbb{T}^s, \mathbb{W}^{s+1})$,
	$V_{\bm{r},\bp} \in \mathcal{L}(\mathbb{U}^s, \mathbb{Y}^{s+1})$
	define compact operators.
\end{enumerate}
\end{lemma}
\begin{proof} 
For the sake of brevity, we assume that $\cP$ is scalar and consider only spaces $T^s,W^{s+1}$, as either vector $\cP$ or the case of spaces $U^s,Y^{s+1}$
follow verbatim.

We start by proving item (i), i.e.~when $\bp = \bm{r}$ for $V_{\bm{r},\bm{r}}$. To this end, let us recall the decomposition of the fundamental solution \eqref{eq:kernelsplit} and define $V_{\bm{r},\bm{r}}^R$
(resp.~$V_{\bm{r},\bm{r}}^S$) for the integral operator with kernel
$G_{\bm{r}}^R$ (resp.~$G_{\bm{r}}^S$). Hence, we first proceed to show that $V_{\bm{r},\bm{r}}^R$ fulfills the assumptions of Theorem \ref{thrm:abstractholm}.
\begin{itemize}
	\item[(i)]
	The operator $V_{\bm{r},\bm{r}}^R$ satisfies Theorem \ref{thrm:abstractholm}
	with $S \equiv 1$ and $p_k = G_{\bm{r}}^R$.
	\item[(ii)]
	Thus, it follows from Lemma \ref{lemma:selfkernel}
	that there exists $\delta>0$ such that
	$$G_{\bm{r}}^R \in \cmspaceh{m-1}{\alpha}{\iinterv}{\IC},$$
	for each $\bm{r} \in K_\delta$.
	Furthermore, Lemma \ref{lemma:selfkernel} ensures that the map
	$$K \ni \bm{r} \mapsto G_{\bm{r}}^R \in \cmspaceh{m-1}{\alpha}{\iinterv}{\IC}$$%
	admits a bounded holomorphic extension into $K_\delta$.
	\item[(iii)]
	By assuming that Condition \ref{condition:smh} holds for a triple $(m,\alpha,s)$,
	it follows from Corollary \ref{corollary:rfop} that for each $\bm{r} \in K_\delta$ one has
	$V_{\bm{r},\bm{r}}^R \in \mathcal{L}\left( T^s,W^{s+1} \right)$, and that, furthermore, it defines
	a compact operator, satisfying
	\begin{align}
	\label{eq:VRbound}
		\norm{V_{\bm{r},\bm{r}}^R }_{\mathcal{L}\left( T^s,W^{s+1}\right)} 
		\lesssim 
		\norm{G^R_{\bm{r}}}_{\cmspaceh{m-1}{\alpha}{\iinterv}{\IC}},
	\end{align}
	where the implied constant is independent of the parametrization
	$\bm{r}: (-1,1) \rightarrow\IR^2$.
	By \eqref{eq:VRbound} and, again, Lemma \ref{lemma:selfkernel}
	the quantity $\|V_{\bm{r},\bm{r}}^R \|_{\mathcal{L}\left( T^s,W^{s+1}\right)}$ 
	is uniformly bounded over $\bm{r} \in K_\delta$. 
\end{itemize}
It follows from Theorem \ref{thrm:abstractholm} that the map
$K \ni\bm{r} \mapsto V^R_{\bm{r},\bm{r}} \in \mathcal{L}(T^s, \W^{s+1})$ admits a bounded 
holomorphic extension onto $K_\delta$.

Now, let us consider $G_{\bm{r}}^S$ and decompose it as follows:
$G_{\bm{r}}^S= G_{\bm{r}}^{S,1}+G_{\bm{r}}^{S,2}$ with:
\begin{align}
	G_{\bm{r}}^{S,1}(t,\tau) & \coloneqq 2F_1(d_{\bm{r}}^2(t,t)) \log|t-\tau|,\\
	G_{\bm{r}}^{S,2}(t,\tau) & \coloneqq 2(F_1(d_{\bm{r}}^2(t,\tau))-F_1(d_{\bm{r}}^2(t,t)))\log|t-\tau|.
\end{align}
The integral operator with kernel $G_{\bm{r}}^{S,1}$
(resp. $G_{\bm{r}}^{S,2}$) is denoted by $V^{S,1}_{\bm{r},\bm{r}}$
(resp. $V^{S,2}_{\bm{r},\bm{r}}$). 
Observe that for each $\bm{r}$ it holds 
$d^2_{\bm{r}}(t,t)=0$ for all $t\in (-1,1)$, and thus, 
$F_1(d^2_{\bm{r}}(t,t)) = F_1(0)$.
Consequently, $V^{S,1}_{\bm{r},\bm{r}}$ is
independent of the parametrization and
it follows from Corollary \ref{corollary:lfop} that
$V^{S,1}_{\bm{r},\bm{r}} \in \mathcal{L}(T^s,W^{s+1})$.
The map $\bm{r} \mapsto V^{S,1}_{\bm{r},\bm{r}}$ is trivially
holomorphic, as it is constant with respect to the parametrization $\bm{r}$.
Moreover, in the representation of the fundamental solution we
have assumed that $F_1(0) \neq 0 $. Hence, $V_{\bm{r},\bm{r}}^{S,1}$ is
invertible from $T^s$ into $W^{s+1}$ (\emph{cf}.~\cite{JEREZHANCKES2011547}).
Notice that when $\cP$ is a vector-valued operator the conclusion still holds as $F_1(0)$
is assumed to be an invertible matrix. 
However, for the pair $U^s,Y^{s+1}$ this does not hold as
the associated integral operator to the logarithmic term is not an invertible operator on the mentioned spaces.

We proceed to show that the operator $V^{S,2}_{\bm{r},\bm{r}}$
fulfils the assumptions of Theorem \ref{thrm:abstractholm}.
\begin{itemize}
	\item[(i)]
	Taylor's theorem yields
	\begin{equation}
		F_1(d_{\bm{r}}(t,\tau)^2)-F_1(0)
		= 
		d_{\bm{r}}(t,\tau)^2
		\int_{0}^1 
		F'_1(\eta d^2_{\bm{r}}(t,\tau))
		\text{d}\eta,
	\end{equation}
	where we have used $d_{\bm{r}}(t,t) = 0$
	for all $t\in (-1,1)$.
	By expanding the distance function,
	for each $(t,\tau) \in (-1,1) \times (-1,1)$ we obtain
	\begin{align}
		F_1(d_{\bm{r}}(t,\tau)^2)-F_1(0)
		= 
		(t-\tau)^2 
		f_{\bm{r}}(t,\tau)
	\end{align}
	where
	\begin{align}
 \begin{split}
		f_{\bm{r}}(t,\tau)
		\coloneqq
		\left(\int_{0}^1 \bm{r}'(t+\eta(\tau-t))\text{d}\eta \right)
		\cdot 
		\left(\int_{0}^1 \bm{r}'(t+\eta(\tau-t))\text{d}\eta \right) \\ \times
		\int_{0}^1 
		F'_1(\eta d^2_{\bm{r}}(t,\tau))
		\text{d}\eta.
  \end{split}
	\end{align}
	We obtain the following representation 
	of $G_{\bm{r}}^{S,2}$
	$$
		G_{\bm{r}}^{S,2}(t,\tau) 
		= 
		(t-\tau)^2 \log|t-\tau| f_{\bm{r}}(t,\tau).
	$$
	Consequently, for each $\bm{r} \in K$, the operator $V^{S,2}_{\bm{r},\bm{r}}$
	fits the framework of Theorem \ref{thrm:abstractholm}
	with $S(t) = t^2\log|t|$ and $p_{k}(t,\tau) = f_{\bm{r}}(t,\tau)$.
	\item[(ii)]
	It follows from Lemma \ref{lemma:selfkernel} that the map
	\begin{equation}
		K\ni\bm{r} 
		\mapsto 
		f_{\bm{r}} \in \cmspaceh{m-1}{\alpha}{\iinterv}{\IC}
	\end{equation}
	admits a bounded holomorphic extension into $K_\delta$
	for some $\delta>0$.
	\item[(iii)]
	Assume that the triple $(m,\alpha,s)$ satisfies
	Condition \ref{condition:smh}. 
	Then, by
	Corollary \ref{corollary:sfop}, for each $\bm{r} \in K$
	we have that $V_{\bm{r},\bm{r}}^{S,2} \in \mathcal{L}\left(T^s,W^{s+1}\right)$, and,
	furthermore, $V_{\bm{r},\bm{r}}^{S,2}$ 
	defines a compact operator satisfying the bound
	\begin{equation}\label{eq:bound_V_S_2}
		\norm{V_{\bm{r},\bm{r}}^{S,2}}_{\mathcal{L}\left(T^s,W^{s+1}\right)} 
		\lesssim 
		\norm{f_{\bm{r}}}_{\cmspaceh{m-1}{\alpha}{\iinterv}{\IC}}.
	\end{equation}
	The right-hand side of \eqref{eq:bound_V_S_2}
	is uniformly bounded on $K_\delta$ as a consequence
	of Lemma \ref{lemma:selfkernel}.
\end{itemize}
It follows from Theorem \ref{thrm:abstractholm} that the map
$K \ni\bm{r} \mapsto V^{S,2}_{\bm{r},\bm{r}} \in \mathcal{L}(T^s, W^{s+1})$ admits a bounded 
holomorphic extension onto $K_\delta$.

Lastly, we obtain the holomorphic extension of $V_{\bm{r},\bm{r}}$ by acknowledging that 
\begin{equation}
	V_{\bm{r},\bm{r}} 
	=
	V_{\bm{r},\bm{r}}^{R}+V_{\bm{r},\bm{r}}^{S,1} + V_{\bm{r},\bm{r}}^{S,2},
\end{equation}
since the three operators on the right-hand side have holomorphic extension at least in $K_\xi$, for $0<\xi < \delta$. On the other hand $V_{\bm{r},\bm{r}}^{R}, V_{\bm{r},\bm{r}}^{S,2}$ are compact operators, and since $V_{\bm{r},\bm{r}}^{S,1}$ is invertible,
$V_{\bm{r},\bm{r}}$ is Fredholm of index zero.

The proof of the second part of the lemma is proved as with the part involving
$V^R_{\bm{r},\bm{r}}$, but by using Lemma \ref{lemma:crosskernel} instead
of Lemma \ref{lemma:selfkernel}. We skip it for the sake of brevity.
\end{proof}

We also consider the generic hyper-singular operator interaction
and establish the sought shape holomorphy property for this type
of BIOs.
To do so, we employ the previously assumed Maue's representation formula \eqref{eq:mauesrep}, which takes the form
\begin{align}
	\left(
		W_{\bm{r},\bp} u
	\right) (t)
	= 
	\frac{\text{d}}{\text{d}t}
	\int_{-1}^1 
	G(\bm{r}(t),\bp(\tau))\frac{\text{d}}{\text{d}\tau}u(\tau) 
	\text{d}\tau
	+ 
	\int_{-1}^1
	\widetilde{G}(\bm{r}(t),\bp(\tau))u(\tau)
	\text{d}\tau
\end{align}
where again we have
$W_{\bm{r}_i,\bm{r}_j} = (W_{\bm{r}_1,\hdots,\bm{r}_M})_{i,j}$, 
and also the next result.

\begin{lemma}
\label{lemma:wcomp} 
Assume that Condition \ref{condition:smh} holds for some
$m \in \IN$, $\alpha \in [0,1]$, and $s \in \IR$.
\begin{enumerate}
	\item[(i)] 
	Let $K$ be an $(m,\alpha)$-admissible set
	of arc parametrizations. Then, there exists $\delta>0$, depending only on $K$
	and satisfying Condition \ref{condition:Qself}, such that
	\begin{align}
		K\ni\bm{r}
		\mapsto W_{\bm{r},\bm{r}} 
		\in 
		\mathcal{L}(\mathbb{U}^s, \mathbb{Y}^{s-1}),
	\end{align}
	admits a bounded holomorphic extensions into $K_\delta$.
	Furthermore, for any $\bm{r} \in K_{\delta}$ one has that
	$W_{\bm{r},\bm{r}} \in \mathcal{L}(\mathbb{U}^s, \mathbb{Y}^{s-1})$
	is a Fredholm operator of index zero.
	\item[(ii)]
	Let $K^1,K^2$ be two $(m,\alpha)$-admissible 
	sets of arc parametrizations satisfying
	Assumption \ref{assump:admissible_arc_param}.
	Then, there exist $\delta_1, \delta_2>0$,
	satisfying Condition \ref{condition:dcross},
	such that the map
	\begin{align}
			K^1 \times K^2 
		\ni
		(\bm{r},\bp)
				\mapsto W_{\bm{r},\bp} 
		\in 
		\mathcal{L}(\mathbb{U}^s, \mathbb{Y}^{s-1}),
	\end{align}
	admits a bounded holomorphic extensions into 
	$K_{\delta_1} \times K_{\delta_2}$.
	Moreover, 
	for any $(\bm{r},\bp) \in K^1_{\delta_1} \times K^2_{\delta_2}$,
	the map $W_{\bm{r},\bp} \in\mathcal{L}(\IU^s, \IY^{s-1})$ is compact. 
\end{enumerate}
\end{lemma}
\begin{proof}
Again, we restrict ourselves to the scalar case. Except for the Fredholm order, the proof follows directly from Maue's representation formula, the mapping properties of the derivative operators \eqref{eq:devprop}---also independent of the parametrizations--- and the arguments of
Lemma \ref{lemma:vcomp}. 

To show the Fredholm order, we use the decomposition of the hyper-singular operator $W_{\bm{r},\bp}=W^1_{\bm{r},\bp}+W^2_{\bm{r},\bp}$, with 
\begin{align}
	\left(W^1_{\bm{r},\bp}u\right)(t)
	&
	\coloneqq
	\frac{\text{d}}{\text{d}t}
	\int_{-1}^1 
	G(\bm{r}(t),\bp(\tau))\frac{\text{d}}{\text{d}\tau}u(\tau) 
	\text{d}\tau, \\ 
	\left(W^2_{\bm{r},\bp}u\right)(t)
	&
	\coloneqq
	\int_{-1}^1
	\widetilde{G}(\bm{r}(t),\bp(\tau))u(\tau)
	\text{d}\tau.
\end{align}
For $W^1_{\bm{r},\bp}$, we argue as in Lemma \ref{lemma:vcomp}.
We decompose this operators into three parts: 
two of them are compact for the same reasons of the
previous lemma, and the remaining part is 
\begin{align}\label{eq:w3}
 	2F_1(0) \frac{\text{d}}{\text{d}t} 
	\int_{-1}^1 \log|t-s| \frac{\text{d}}{\text{d}s}
	u(\tau)
	\text{d}s,   
\end{align}
The factor $2F_1(0)$ 
is assumed to be invertible, and the integral operator
is the standard hyper-singular one for the Laplace equation. 
Hence, the operator in \eqref{eq:w3} is invertible
as a map in $\mathcal{L}(\IU^s,\IY^{s-1})$ 
(\emph{cf}.~\cite{JEREZHANCKES2011547}).

The operator $W^2_{\bm{r},\bp}$ can be analyzed as in
Lemma \ref{lemma:vcomp}, so as to find that $W^2_{\bm{r},\bp} \in \mathcal{L}\left( \IU^s,\IY^{s+1}\right)$.
Therefore, by the compact embedding of the corresponding spaces (see Section \ref{sec:FunctionalSpaces}) we have that $W^2_{\bm{r},\bp} \in \mathcal{L}\left( \IU^s,\IY^{s-1}\right)$ is a compact operator.
\end{proof}

\begin{remark}
\label{remark:evenlessreg}
In practice, the term $\widetilde{G}$ of Maue's representation formula includes a factor involving normal vectors. This ensures that the corresponding functions $F_1,F_2$ have the structure described in Remark \ref{remark:lessregularitykernel}. Consequently, none of the results needs to be modified.

We can further generalize the structure of the functions $G_1,G_2$ in Remark \ref{remark:lessregularitykernel} by considering the form:
\begin{equation}
G_j(t,\tau) = f_j\left((\bm{r}'(t),\bm{r}''(t),\hdots,\bm{r}^{(n)}(t),\bp'(\tau),\bp''(\tau),\hdots,\bp^{(n)}(\tau)\right)F_j(d_{\bm{r},\bp}^2(t-\tau)),
\end{equation}
where $f_j$ is entire in each coordinate. With the above representation, Condition \ref{condition:smh} is changed to:
\begin{enumerate}
\item[(i)] 
$s > -\frac{1}{2}$, and $s+
\frac{3}{2}+n<m+\alpha$,
\item[(ii)] 
$ s \leq \frac{-1}{2}$, and $\frac{1}{2}+n -s < m+\alpha$.
\end{enumerate}  
\end{remark}
\subsection{Multiple arcs ($M>1$)}
\label{sec:multiple_arcs}
Lemmas \ref{lemma:vcomp} and \ref{lemma:wcomp} ensure that for every pair of arcs $\bm{r}, \bp$ there exists a region---depending on the arcs---such that the weakly- and hyper-singular BIOs have a holomorphic extension. Now we return to the original problem (Sections \ref{sec:bvproblem}, and \ref{sec:bif}) to prove the BIOs' holomorphic extension for the interaction of $M >1$ arcs (\emph{cf.}~Theorem \ref{thrm:biohlmextension}). With this, we obtain the holomorphic extension of the so called domain to solution map for the problem presented in Section \ref{sec:bvproblem}.


\begin{condition}
\label{condition:deltaOps}
Consider $M$ different $(m,\alpha)$-admissible sets of parametrizations $K^1,\hdots,K^M$ satisfying 
Assumption \ref{assump:admissible_arc_param}.
Let $\delta_1,\hdots,\delta_M$ be $M$ strictly positive real numbers,
satisfying:
\begin{enumerate}
	\item[(i)] 
	Each $\delta_j$ is such that
	Condition \ref{condition:Qself} is
	fulfilled in the compact set $K^j$, for $j=1,\hdots,M$.
	\item[(ii)]
	For each $(\delta_i, \delta_j)$
	with $i,j \in \{1,\hdots,M\}$ and $i \neq j$,
	Condition \ref{condition:dcross} is fulfilled
	in $K^i\times K^j$.
\end{enumerate} 
\end{condition}



\begin{theorem}
\label{thrm:biohlmextension}
Let $s \in \IR$, $m \in \IN$ and $\alpha \in [0,1]$
be such that Condition \ref{condition:smh} is fulfilled.
 Let $K^1,\hdots,K^M$ be $M$ $(m,\alpha)$-admissible
sets of parametrizations satisfying 
Assumption \ref{assump:admissible_arc_param}.
Then there exist $\delta_1,\dots,\delta_M>0$ satisfying
Condition \ref{condition:deltaOps} such that
\begin{align}
	K^1 \times \hdots \times K^M 
	\ni
	(\bm{r}_1,\hdots,\bm{r}_M) 
	\mapsto
	\mathbf{V}_{\bm{r}_1, \dots,\bm{r}_M}
	\in \mathcal{L} 
	\left(
		\prod_{j=1}^M \mathbb{T}^s, \prod_{j=1}^M \mathbb{W}^{s+1}
	\right),
	\\
	K^1 \times \hdots \times K^M 
	\ni
	(\bm{r}_1,\hdots,\bm{r}_M)
	\mapsto
	\mathbf{W}_{\bm{r}_1, \dots,\bm{r}_M}
	\in
	\mathcal{L}
	\left(
		\prod_{j=1}^M \mathbb{U}^s, \prod_{j=1}^M \mathbb{Y}^{s-1}
	\right),
\end{align}
admit bounded holomorphic extensions into
$K^1_{\delta_1}\times\hdots{\times} K^M_{\delta_M}$.
\end{theorem}

\begin{proof}
We prove only the result for the weakly-singular BIO
as the hyper-singular case follows similarly. Our first observation is that one can write
\begin{align}
	\mathbf{V}_{\bm{r}_1,\hdots,\bm{r}_M} 
	=
	\begin{pmatrix}
		V_{\bm{r}_1,\bm{r}_1} &V_{\bm{r}_1,\bm{r}_2}& \hdots &V_{\bm{r}_1,\bm{r}_M}\\
		V_{\bm{r}_2,\bm{r}_1}& V_{\bm{r}_2,\bm{r}_2}& \hdots &V_{\bm{r}_2,\bm{r}_M}\\
		\vdots & \ddots &\hdots & \vdots\\
		V_{\bm{r}_M,\bm{r}_1} & V_{\bm{r}_M,\bm{r}_2} &\hdots & V_{\bm{r}_M,\bm{r}_M}
	\end{pmatrix}.
\end{align}

It follows from Lemma \ref{lemma:vcomp}
that there exist $\delta_1,\dots,\delta_M>0$
satisfying Condition \ref{condition:deltaOps}
such that the maps
\begin{align}
	K^j
	\ni
	\bm{r}_j &
	\mapsto 
	V_{\bm{r}_j,\bm{r}_j}
	\in 
	\mathcal{L}\left(\mathbb{T}^s,\mathbb{W}^{s+1}\right), \\
	K^i \times K^j
	\ni 
	(\bm{r}_i,\bm{r}_j) &
	\mapsto 
	V_{\bm{r}_i,\bm{r}_j}
	\in 
	\mathcal{L}\left(\mathbb{T}^s,\mathbb{W}^{s+1}\right), \quad i \neq j,
\end{align}
admit bounded holomorphic extension
into $K^j_{\delta_j}$ and $K^i_{\delta_i} \times K^j_{\delta_j}$, respectively. 
Since each component has a holomorphic extension, by defining the norms for 
$\prod_{j=1}^M \mathbb{T}^s$ 
and 
$\prod_{j=1}^M \mathbb{W}^{s+1}$ as the standard Euclidean norm of a Cartesian product space, we directly deduce that 
\begin{align}
	K^1 \times \hdots K^M 
	\ni
	(\bm{r}_1,\hdots,\bm{r}_M) 
	\mapsto 
	\mathbf{V}_{\bm{r}_1,\hdots,\bm{r}_M}
	\in 
	\mathcal{L} 
	\left(
		\prod_{j=1}^M \mathbb{T}^s, \prod_{j=1}^M \mathbb{W}^{s+1}
	\right)
\end{align}
admits a bounded holomorphic extension into
$K^1_{\delta_1} \cdots \times K^M_{\delta_M}$.
\end{proof}

From this last result, by Theorem \ref{thrm:abtractinverse},
and assuming that the right-hand sides of the BVPs are given by entire functions  (Section \ref{sec:bvproblem}),
we conclude that $\bla_{\bm{r}_1,\hdots,\bm{r}_M}$ and $\bmu_{\bm{r}_1,\hdots,\bm{r}_M}$,
solution to the Dirichlet and Neumann problems, respectively, depend
holomorphically upon perturbations of arc parametrizations
$\bm{r}_1,\hdots,\bm{r}_M$.
 
\begin{theorem}
\label{thrm:solutionholomrext}
Under the same hypothesis of Theorem \ref{thrm:biohlmextension}, 
there exists $\eta>0$ such that the maps
\begin{align}
\label{eq:d2s_map_dir}
	K_1 \times \hdots \times K_M \ni
	(\bm{r}_1,\hdots,\bm{r}_M) 
	\mapsto  
	\bla_{\bm{r}_1,\hdots,\bm{r}_M}
	\in \prod_{j=1}^M
	\mathbb{T}^s
\end{align}
and
\begin{align}
\label{eq:d2s_map_neu}
	K_1 \times \hdots \times K_M 
	\ni
	(\bm{r}_1,\hdots,\bm{r}_M)
	\mapsto  
	\bmu_{\bm{r}_1,\hdots,\bm{r}_M} \in \prod_{j=1}^M 
	\mathbb{U}^s
\end{align}
admit bounded holomorphic extensions into
$K_{\eta}^1 \times \hdots K_{\eta}^M$,
where for each 
$(\bm{r}_1,\hdots,\bm{r}_M)  \in K^1 \times \hdots \times K^M$
we have that
$\bla_{\bm{r}_1,\hdots,\bm{r}_M}$ and $\bmu_{\bm{r}_1,\hdots,\bm{r}_M}$
are the boundary solutions of the Dirichlet and Neumman
problems stated in \eqref{eq:bios}.
\end{theorem}

\begin{proof}
As before, we prove only the result for the weakly-singular BIO, though remarks are given whenever the proof differs for the
hyper-singular BIO. 

The proof relies on Theorem \ref{thrm:abtractinverse}. Firstly, we need to ensure that for each 
$(\bm{r}_1,\hdots,\bm{r}_M)  \in K^1 \times \hdots K^M$
the maps introduced in \eqref{eq:d2s_map_dir}
and \eqref{eq:d2s_map_neu} are well defined. To this end, we use the block-wise decomposition of the
weakly-singular BIO defined on multiple disjoint arcs into diagonal and
off-diagonal components, i.e. 
\begin{align}\label{eq:block_decom_V}
	\mathbf{V}_{\bm{r}_1,\hdots,\bm{r}_M}
	=
	\begin{pmatrix}
		V_{\bm{r}_1,\bm{r}_1} & \hdots &0\\
		\vdots& \ddots& \vdots\\
		0 &\hdots & V_{\bm{r}_M, \bm{r}_M}
	\end{pmatrix}
	+
	\begin{pmatrix}
		0 & \hdots &V_{\bm{r}_1,\bm{r}_M}\\
		\vdots& \ddots& \vdots\\
		V_{\bm{r}_M,\bm{r}_1} &\hdots & 0
	\end{pmatrix}
\end{align}
An equivalent decomposition can be stated
for the hyper-singular BIO. It follows from Lemma \ref{lemma:vcomp}---Lemma \ref{lemma:wcomp} for the hyper-singular BIO---
that the diagonal part, i.e. the first summand in \eqref{eq:block_decom_V},
is composed of index zero Fredholm operators, while the off-diagonal 
include compact operators.
Recalling the definition of a Cartesian product space,
the block operators are Fredholm of index zero. 
From the Fredholm property, we obtain that 
$\mathbf{V}_{\bm{r}_1,\hdots,\bm{r}_M} \in \mathcal{L} \left( \prod_{j=1}^M \mathbb{T}^{s}, \prod_{j=1}^M \mathbb{W}^{s+1} \right)$ 
---also
$\mathbf{W}_{\bm{r}_1,\hdots,\bm{r}_M} \in \mathcal{L} \left( \prod_{j=1}^M U^{s}, \prod_{j=1}^M Y^{s-1} \right)$---is invertible provided that is injective.
As mentioned, this is equivalent to the unisolvence of the 
corresponding volume problem presented in Section \ref{sec:bvproblem} (\emph{cf}.~\cite{JHP20}).

Hence, for each
$(\bm{r}_1,\hdots, \bm{r}_M) \in K^1\times \cdots \times K^M$,
the operators $\mathbf{V}_{\bm{r}_1,\hdots,\bm{r}_M}$,
$\mathbf{W}_{\bm{r}_1,\hdots,\bm{r}_M}$ are invertible mapping as follows
\begin{align*}
    \begin{split}
({\mathbf{V}_{\bm{r}_1,\hdots,\bm{r}_M}})^{-1} 
	\in 
	\mathcal{L} 
	\left(
		\prod_{j=1}^M \mathbb{W}^{s+1}, \prod_{j=1}^M \mathbb{T}^s 
	\right), \\
(\mathbf{W}_{\bm{r}_1,\hdots,\bm{r}_M})^{-1} 
	\in 
	\mathcal{L} 
	\left( 
		\prod_{j=1}^M \mathbb{Y}^{s-1}, \prod_{j=1}^M \mathbb{U}^s 
	\right).        
    \end{split}
\end{align*}
Thus, we have proved that the maps introduced in \eqref{eq:d2s_map_dir}
and \eqref{eq:d2s_map_neu} are well-defined. It follows straightforwardly from Theorems \ref{thrm:abtractinverse}
and \ref{thrm:biohlmextension} the there exists $\eta>0$ such the maps
\begin{align}	
	K^1 \times \cdots \times K^M 
	\ni
	(\bm{r}_1,\hdots,\bm{r}_M) 
	&
	\mapsto 
	(\mathbf{V}_{\bm{r}_1,\dots,\bm{r}_M})^{-1} 
	\in 
	\mathcal{L} 
	\left(
		\prod_{j=1}^M \mathbb{W}^{s+1}, \prod_{j=1}^M \mathbb{T}^s
	\right)\\
	K^1 \times \cdots \times K^M 
	\ni
	(\bm{r}_1,\hdots,\bm{r}_M) 
	&
	\mapsto
	\left(
		\mathbf{W}_{\bm{r}_1,\dots,\bm{r}_M}
	\right)^{-1}
	\in
	\mathcal{L} 
	\left(
		\prod_{j=1}^M \mathbb{Y}^{s-1}, \prod_{j=1}^M \mathbb{U}^s 
	\right)
\end{align}
admit bounded holomorphic extension holomorphic into $K_{\eta}$.
Recall that the right-hand sides of the boundary integral formulations 
considered are of the form:
\begin{align}
	({f}^D_{\bm{r}_1,\hdots,\bm{r}_M})_j 
	= 
	f^D \circ \bm{r}_j
	\quad
	\text{and} 
	\quad
	({f}^N_{\bm{r}_1,\hdots,\bm{r}_M})_j 
	= 
	f^N \circ \bm{r}_j,
\end{align}
and we assumed that $f^D$, $f^N$ are entire functions.
From this representation, arguing as in Lemma \ref{lemma:cmdecay2},
but for univariate functions, one can check that
Condition \ref{condition:smh} ensures that 
\begin{align}
	{f}^D_{\bm{r}_1,\hdots,\bm{r}_M} 
	\in 
	\prod_{j=1}^M \mathbb{W}^{s+1}, 
	\quad {f}^N_{\bm{r}_1,\hdots,\bm{r}_M} 
	\in 
	\prod_{j=1}^M \mathbb{Y}^{s+1} 
	\subset 
	\prod_{j=1}^M \mathbb{Y}^{s-1}.
\end{align}
Furthermore, once again using that $f^D$, $f^N$ are entire functions, it is direct to see that the maps $(\bm{r}_1,\hdots,\bm{r}_M) \mapsto \mathbf{f}^D_{\bm{r}_1,\hdots,\bm{r}_M}$, and  $(\bm{r}_1,\hdots,\bm{r}_M) \mapsto \mathbf{f}^N_{\bm{r}_1,\hdots,\bm{r}_M}$ 
admit bounded holomorphic extensions to any region. The final result follows by composition of maps with holomorphic extensions. 
\end{proof}

\begin{remark}\label{rmk:linear_functionals}
Theorem \ref{thrm:solutionholomrext} enables us to obtain holomorphic
extensions for some linear functionals. In particular, if we consider linear
functionals of the form  
\begin{equation}
	\left(L\bu\right)(\x) 
	= 
	\int_{-1}^{1}
	\boldsymbol{\vartheta}(\x,\bm{r}_1(t),\hdots,\bm{r}_M(t)) 
	\cdot 
	\bu(t) 
	\emph{\text{d}}t,
	\quad
	\x \in \IR^2,
\end{equation}
where $\bu = \bla $ or $\bu = \bmu$, solutions of the corresponding BIEs (Problem \ref{prob:BIEs}), and $\boldsymbol{\vartheta}$
is entire on each coordinate, except possibly on
the first one. The holomorphic extension
of this functional is proven by showing that it is a composition of holomorphic functions. First, we notice that by Theorem \ref{thrm:abstractholm}, $(\bm{r}_1, \hdots \bm{r}_M) \rightarrow L \cdot$, has a bounded holomorphic extension, and secondly, by the previous theorem,, $\bu$ also have a bounded holomorphic extension. Hence, $L u$  have bounded  holomorphic extension. 
\end{remark}

\subsection{Parametric Holomorphy of the Domain-to-Solution Map}
\label{sec:parametrichlmrf}

Throughout this section we denote by 
$\bm{r}^0_1,\bm{r}^1_2,\hdots, \bm{r}^0_M$ 
a collection of $M$ arc parametrizations,
each of them contained in $\rgeoh{m}{\alpha}$, and such that no crossing among them occurs. Specifically, we consider the next affine-parametric arc 
parametrizations:
\begin{align}\label{eq:affine_parametric_arcs_basis}
	\bm{r}_{j,\by_j}
	= 
	\bm{r}^0_j 
	+ 
	\sum_{n=1}^\infty y^{n}_j
	\bm{r}^n_j,
	\quad
	j = 1,\dots,M,
	\quad
	\by_j \coloneqq (y^n_j)_{n\in \IN}\in \sU,
\end{align}
where, for each $j\in \{1,\dots,M \}$,
the sequence 
$\{\bm{r}^n_j \}_{n\in \IN}\subset \cmspaceh{m}{\alpha}{(-1,1)}{\IR^2}$.
For each $n\in \IN$ and $j=1,\dots,M$, 
let us set
\begin{equation}\label{eq:definition_b_j}
	b^n_j\coloneqq\norm{\bm{r}^n_j}_{\cmspaceh{m}{\alpha}{(-1,1)}{\IR^2}}
	\quad
	\text{and}
	\quad
	\boldsymbol{b}_{j}= \{b^n_j\}_{n\in \IN}.
\end{equation}
In order to restrict ourselves to admissible
geometric configurations, i.e.~ arc parametrization
satisfying $\bm{r}_{j,\by_j} \in \rgeoh{m}{{\alpha}}$ for $j=1,\hdots,M$
and each $\by \in \sU$,
such that no crossings or intersections occur as well as to adopt the framework of Section \ref{section:ShapeHolormsGeneral},
we work under the following assumptions. 

\begin{assumption}
\label{assump:geoparam}
We assume that
\begin{enumerate}
	\item[(i)] 
	For each $j\in \IN$
	the sequence
	$\boldsymbol{b}_{j} \in \ell^p(\IN)$ for some $p\in (0,1)$.
	\item[(ii)] 
	There exists a single $\zeta\in(0,1)$
	such that, for each $j \in \{1,\hdots,M\}$, it holds that
	\begin{align}
		\sup_{t\in (-1,1)} 
		\sum_{n=1}^\infty 
		\|(\bm{r}_j^n)'(t)\| 
		\leq
		\zeta
		\inf_{t \in (-1,1) }
		\|(r^0_j)'(t)\| .
	\end{align}
\item[(iii)]
There exists $\eta \in (0,1)$
such that, for any $i,j \in \{1,\dots,M\}$
and for each $\by_i,\by_j \in \sU$, one has
\begin{align}
	\norm{
		\sum_{n=1}^\infty 
		y^n_i \bm{r}^n_i(t)
		-
		y^n_j \bm{r}^n_j(\tau)
	}
	\leq
	\eta
	\inf_{(t,\tau) \in (-1,1)\times (-1,1)}
	\norm{
		\bm{r}^0_i(t)
		-
		\bm{r}^0_j(\tau)
	}.
\end{align}
\end{enumerate}
\end{assumption}


In the following, for $j\in \{1,\dots,M\}$ we set
\begin{align}\label{eq:K_j_param}
	K_j
	\coloneqq
	\left\lbrace 
		\bm{r}_{j,\by}
		\in
		\cmspaceh{m}{\alpha}{(-1,1)}{\IR^2}: 
		\bm{r}_{j,\by} 
		= 
		\bm{r}^0_j 
		+ 
		\sum_{n=1}^\infty 
		y^n \bm{r}^n_j, 
		\ 
		\by = \{y^n\}_{n \in \IN} 
		\in 
		\normalfont\text{U}
	\right
	\rbrace.
\end{align} 
Observe that due to item (i) in Assumption 
\ref{assump:geoparam} the series in 
\eqref{eq:affine_parametric_arcs_basis} converges
absolutely and uniformly with respect to $\by\in \sU$. Moreover, one can now obtain a proper set of arc parametrizations as shown below.

\begin{lemma}
Let $K_1,\dots,K_M \subset \cmspaceh{m}{\alpha}{(-1,1)}{\IR^2}$
be as in \eqref{eq:K_j_param}
for some $m \in \IN$ and $\alpha \in [0,1]$,
and let Assumption \ref{assump:geoparam} be satisfied.
Then $K_1,\dots,K_M$ are $(m,\alpha)$-admissible
arc parametrizations in the sense of Definition
\ref{assump:admissible_arc_param} satisfying 
Assumption \ref{assump:admissible_arc_param}.
\end{lemma}%
\begin{proof}

We start by proving that for each $\by \in \sU$ the arc parametrization 
$\bm{r}_{j,\by}: (-1,1) \rightarrow \IR^2$ defined as in 
\eqref{eq:affine_parametric_arcs_basis} renders an element
of $\rgeoh{m}{\alpha}$. 
The proof follows substantially
that of \cite[Lemma 6.2]{henriquez2021_thesis}, however
we include the details for the sake of completeness. 

For $t,\tau \in (-1,1)$, we directly have that 
$$
 \norm{\bm{r}_{j,\by}(t) -\bm{r}_{j,\by}(\tau)} = 
\norm{\bm{r}_j^0(t) - \bm{r}_j^0(\tau) + \sum_{n=1}^\infty y^n(\bm{r}_j^n(t)-\bm{r}_j^n(\tau))}.
$$
Hence, by the Taylor expansion of $\bm{r}_{j,\by}$, for every $j = 0,\hdots M$, there exist $\xi \in (-1,1)$ such that 
\begin{align*}
    \begin{split}
\norm{\bm{r}_{j,\by}(t) -\bm{r}_{j,\by}(\tau)} = (t-\tau)&
\norm{(\bm{r}_j^0)'(\xi) + \sum_{n=1}^\infty y^n(\bm{r}_j^n)'(\xi)}\\
&\geq (t-\tau) \left(\norm{(\bm{r}_j^0)'(\xi)} - \norm{\sum_{n=1}^\infty y^n(\bm{r}_j^n)'(\xi)} \right).
\end{split}
\end{align*}
By Assumption \ref{assump:geoparam} (ii), we obtain
$$
\norm{(\bm{r}_j^0)'(\xi)} - \norm{\sum_{n=1}^\infty y^n(\bm{r}_j^n)'(\xi)} \geq (1-\zeta)\inf_{t \in (-1,1) }
		\|(r^0_j)'(t)\|>0.
$$
We deduce that $\bm{r}_{j,\by}$ is injective for every $\by \in \sU$, thus it has a global inverse. Furthermore, if we make the same analysis for the tangent vector we have that 
$$
\norm{(\bm{r}_{j,\by})'(t)} \geq 
\left(\norm{(\bm{r}_j^0)'(t)} - \norm{\sum_{n=1}^\infty y^n(\bm{r}_j^n)'(t)} \right) \geq  (1-\zeta)\inf_{t \in (-1,1) }
		\|(r^0_j)'(t)\|>0.
$$
Thus, the tangent vector is nowhere null, and we have that $K_j \subset \rgeoh{m}{\alpha}$. As explained in Section \ref{section:ShapeHolormsGeneral},
it follows from \cite[Lemma 2.7]{CD15} that each $K_i$ is 
a compact subset of $\cmspaceh{m}{\alpha}{(-1,1)}{\IR^2}$,
thus rendering each of those sets $(m,\alpha)$-admissible
arc parametrizations. 

To conclude, we verify that under 
Assumption \ref{assump:geoparam} the sets 
$K_1,\dots,K_M$ fulfill Assumption \ref{assump:admissible_arc_param}.
Using item (iii) in Assumption \ref{assump:geoparam}
For any $i,j \in \{1,\dots,M\}$ and for each $\by_i, \by_j \in \sU$, one has that
\begin{align}
	\norm{
		\bm{r}_{i,\by_i}(t)
		-
		\bm{r}_{j,\by_j}(\tau)
	}
	&
	=
	\norm{
		\bm{r}^0_i(t)
		+ 
		\sum_{n=1}^\infty 
		y_n \bm{r}^n_i(t)
		-
		\bm{r}^0_j(\tau)
		-
		\sum_{n=1}^\infty 
		y_n \bm{r}^n_j(\tau)
	} \\
	&
	\geq
	\snorm{
		\norm{
			\bm{r}^0_i(t)
			-
			\bm{r}^0_j(\tau)
		}
		-
		\norm{
			\sum_{n=1}^\infty 
			y^n_i \bm{r}^n_i(t)
			-
			y^n_j \bm{r}^n_j(\tau)
		}
	}
	\geq
	(1-\eta)>0,
\end{align}
as stated.
\end{proof}
Let us set for each $\by \in \sU$
\begin{align}
\label{eq:paramsol}
	\bla_{\by}
	\coloneqq
	\bla_{\bm{r}_{1,\by_1},\dots,\bm{r}_{M,\by_M}}
	\quad
	\text{and}
	\quad
	\bmu_{\by}
	\coloneqq
	\bmu_{\bm{r}_{1,\by_1},\dots,\bm{r}_{M,\by_M}},
\end{align}
where the elements $\by_1, \by_2, \hdots \by_M \in \sU$ are defined as
$$
(y_j)_n=  y_{j+nM}, \quad j \in \{1,\hdots,M\}, \ n \in \IN.
$$
We will make use of the set $\boldsymbol{K} := K_1 \times\hdots\times K_M \subset \prod_{j=1}^M \cmspaceh{m}{\alpha}{(-1,1)}{\IR^2}$, which can be written as 
\begin{align}
\begin{split}
    \label{eq:Kset}
    \boldsymbol{K} = \left\lbrace \mathbf{k}_{\by} = 
\begin{pmatrix}
    \bm{r}_{1,\by_1}\\
    \bm{r}_{2,\by_2}\\
    \vdots \\
    \bm{r}_{M,\by_M}
\end{pmatrix} = 
\begin{pmatrix}
    \bm{r}_{1}^0\\
    \bm{r}_{2}^0\\
    \vdots \\
    \bm{r}_{M}^0
\end{pmatrix}
+
\sum_{n=1}^\infty 
y_n \mathbf{k}_n, \ \by = \{ y_n\}_{n \in \IN} \in \sU
\right\rbrace,
\end{split}
\end{align}
wherein,
$$
\{\mathbf{k}_n\}_{n \in \IN} = \left\lbrace
\begin{pmatrix}
\bm{r}_1^1\\
0\\
\vdots\\
0
\end{pmatrix}, 
\begin{pmatrix}
0\\
\bm{r}_2^1\\
\vdots\\
0
\end{pmatrix}, 
 \hdots 
 \begin{pmatrix}
0\\
0\\
\vdots\\
\bm{r}_M^1
\end{pmatrix},
\begin{pmatrix}
\bm{r}_1^2\\
0\\
\vdots\\
0
\end{pmatrix}, 
\begin{pmatrix}
0\\
\bm{r}^2_2\\
\vdots\\
0
\end{pmatrix}, \hdots 
\begin{pmatrix}
0\\
0\\
\vdots\\
\bm{r}^2_M
\end{pmatrix},\hdots
\right\rbrace
$$
or more rigorously  $
\mathbf{k}_n = \bm{r}_{j(n)}^{\left\lceil \frac{n}{M}\right\rceil}\mathbf{e}_{j(n)}
$, with $\left\lceil \frac{n}{M}\right\rceil$ the upper integer part of $ \frac{n}{M}$,
and $j(n)$ is one plus the residual of the integer division
$\frac{n-1}{M}$.
We also define the parameter norm set:
\begin{align}
\label{eq:bsec}
\boldsymbol{b} :=\{\|\mathbf{k}_n\|_{\prod_{j=1}^M \cmspaceh{m}{\alpha}{(-1,1)}{\IR^2}}\}_{n \in \IN} .
\end{align}
\begin{theorem}
\label{thrm:bpeholomr}
Let Condition \ref{condition:smh} hold for some
$m \in \IN$, $\alpha \in [0,1]$, and $s \in \IR$.
Let Assumption \ref{assump:geoparam} be satisfied 
with $\boldsymbol{b}_{j}$ for $j=1,\dots,M$ as in 
\eqref{eq:definition_b_j} and for some $p\in(0,1)$.
Then the maps
\begin{align}
	\sU
	\ni
	\by
	\mapsto 
	\bla_{\by}
	\in 
	\prod_{j=1}^M 
	\mathbb{T}^s,
	\quad
	\text{and}
	\quad
	\sU
	\ni
	\by
	\mapsto 
	\bmu_{\by}
	\in 
	\prod_{j=1}^M 
	\mathbb{U}^s
\end{align}
are $(\boldsymbol{b},p,\varepsilon)$-holomorphic
for some $\varepsilon>0$, $p\in (0,1)$ as in 
Assumption \ref{assump:geoparam}, and $\boldsymbol{b}$ as in \eqref{eq:bsec}. Also, these maps are
continuous when $\sU$ is equipped 
with the product topology.
\end{theorem}

\begin{proof}
Being a direct consequence of Theorem \ref{theorem:generalbpeholomrfsm}, we only need to verify that the hypotheses are satisfied. The role of the compact set $K$ of Theorem \ref{theorem:generalbpeholomrfsm} is played by $\boldsymbol{K}$, which has the desire form and is compact due to the previous Lemma. Our definition of $\boldsymbol{b}$ \eqref{eq:bsec} coincides with that of Theorem \ref{theorem:generalbpeholomrfsm}, and we also have that $\boldsymbol{b} \in \ell^p(\IN)$ by  Assumption \ref{assump:geoparam} (i). 
Finally, we observe that by Theorem \ref{thrm:solutionholomrext}, the maps
\begin{align*}
\begin{split}    
	\boldsymbol{K}
	\ni
	(\bm{r}_1,\hdots,\bm{r}_M) 
	\mapsto  
	\bla_{\bm{r}_1,\hdots,\bm{r}_M}
	\in \prod_{j=1}^M \mathbb{T}^s, \\
		\boldsymbol{K}
	\ni
	(\bm{r}_1,\hdots,\bm{r}_M)
	\mapsto  
	\bmu_{\bm{r}_1,\hdots,\bm{r}_M} \in \prod_{j=1}^M \mathbb{U}^s
 \end{split}
\end{align*}
have bounded holomorphic extensions. 
Moreover, by definition $\bla_{\by}, \bmu_{\by}$, we have that
 $\bla_{\by} = \bla_{\mathbf{k}_{\by}}$, 
$\bmu_{\by} = \bmu_{\mathbf{k}_{\by}}$, where $\bla_{\mathbf{k}_{\by}}, \bmu_{\mathbf{k}_{\by}}$ are the corresponding elements $f(k_{\by})$ in  Theorem \ref{theorem:generalbpeholomrfsm}. The results then follow from Theorem \ref{theorem:generalbpeholomrfsm} as stated.
\end{proof}

\begin{remark}
The result stated in Theorem \ref{thrm:bpeholomr}
enables us to conclude the parametric holomorphy
of linear functionals acting on $\bla_{\by}$ and $\bmu_{\by}$
(\emph{cr.}~Remark \ref{rmk:linear_functionals}). 
\end{remark}

\section{Applications: Time-Harmonic Acoustic and Elastic Wave Scattering}
\label{sec:examples}
In the following, we consider two particular instances of the operator $\cP$, namely, Helmholtz and elastic wave operators, and check whether the assumptions to guarantee holomorphic extensions of their corresponding BIOs are satisfied.

\subsection{Helmholtz Equation}
\label{sec:helmholtz}
The scalar Helmholtz operator with wavenumber $\kappa \in \IR_+$ is given by $\cP  = -\Delta - \kappa^2$. 
In order to ensure well-posedness, one prescribes the following behavior at infinity
\begin{equation}
	\lim_{\|\x\| \rightarrow \infty}
	\| \x \|^{\frac{1}{2}}\left( \frac{\partial u}{\partial {\|\x\|}}-\imath {\kappa} u\right)
	=
	0,
\end{equation}
known as the Sommerfeld radiation condition.
We refer to \cite{stephan1984augmented,kress1995inverse}
for uniqueness of the Dirichlet problem and to \cite{MONCH1996343} for the Neumann one. For the latter, the co-normal trace becomes the standard Neumann trace, i.e.~$\mathbf{\mathcal{B}}_{\bn}u = \bn \cdot \nabla u  = \partial_{\bn} u$ for $u$ smooth enough, with continuous extensions to Sobolev spaces \cite[Lemma 4.3]{mclean2000strongly}. The fundamental solution of the operator $\cP = -\Delta -k^2$ 
in two-dimensional space is
\begin{align}
\label{eq:Helmholtz}
	G_\kappa(\x,\y) 
	= 
	\frac{i}{4}H^{(1)}_0(\kappa \norm{\x-\y}),\quad {\x,\y\in\mathbb{R}^2,}
\end{align}
where $H^{(1)}_0$ denotes the Hankel function of first kind
and order zero, defined as $H^{(1)}_0(z) = J_0(z) + i Y_0(z)$ for $z\neq 0$, 
where $J_0, Y_0$ are the zeroth order Bessel function of first and second kind, respectively.

\begin{corollary}[Helmholtz case] Consider the Dirichlet and Neumann BIEs (Problem \ref{prob:BIEs}) for the Helmholtz kernel \eqref{eq:Helmholtz} with $\kappa > 0$. Then, the arising BIOs and domain-to-solution maps are shape holomorphic.
\end{corollary}

\begin{proof}
From \cite[9.1.12 and 9.1.13]{abramowitz1965handbook}, one has that 
\begin{equation}
	G_{\kappa}(\x,\y) 
	= 
	-\frac{1}{4\pi}
	\log\|\x-\y\|^2
	J_0(\kappa \|\x -\y\|)
	+ 
	R(\kappa \|\x-\y\|^2),
\end{equation}
where the first kind Bessel function $J_0$, and $R$ are entire functions. Furthermore, from \cite[9.1.12]{abramowitz1965handbook}, one has that $J_0(0) = 1$, and also $J_0(z) = j_0(z^2)$, with $j_0$ being an entire function.  
Hence, the representation of the form of \eqref{eq:funsolgen}  holds with 
$$
	F_1(z) 
	= 
	-\dfrac{1}{2\pi}j_0(\kappa z),
	\quad 
	\text{and}
	\quad 
	F_2(z) 
	= R(\kappa z).
$$
Therefore, shape holomorphy results for the Dirichlet problem follow directly from Theorems \ref{thrm:solutionholomrext} and  \ref{thrm:bpeholomr}. 

For the Neumann case, we make use of the corresponding Maue's representation formula \cite[Corollary 3.3.24]{Sauter:2011}:
\begin{align*}
\begin{split}
	(W_{\Gamma_1,\hdots,\Gamma_M})_{i,j}u(t) 
	=
	& \frac{\text{d}}{\text{d}t} \int_{-1}^1 G(\bm{r}_i(t),\bm{r}_j(\tau)) \frac{\text{d}}{\text{d}\tau} u(\tau) \text{d}\tau \\&- k^2\int_{-1}^1 (\bm{r}_i'(t)\cdot \bm{r}_j'(\tau)) G(\bm{r}_i(t),\bm{r}_j(\tau))u(\tau) \text{d}\tau. 
\end{split}
\end{align*} 
Thus, following the notation of Section \ref{sec:bif} we have that 
$$
\widetilde{G}(\bm{r}_i(t),\bm{r}_j(\tau)) = -k^2 (\bm{r}_i'(t)\cdot \bm{r}_j'(\tau)) G(\bm{r}_i(t),\bm{r}_j(\tau)).
$$
Notice that this function has almost the same structure of $G(\bm{r}_i(t),\bm{r}_j(\tau))$ except for the loss of one degree of regularity because of the factors $\bm{r}_i'$, and $\bm{r}_j'$. However, as mentioned in Remark \ref{remark:lessregularitykernel}, this does not have any impact, and we obtain the corresponding results for the Neumann problem.
\end{proof}

\begin{remark}
For $\kappa=0$ (Laplace operator), the
fundamental solution becomes $G(\x,\y) = -\frac{1}{4\pi} \log \| \x - \by \|^2$. Thus, well-posedness requires a suitable radiation condition. One particular alternative is to impose solutions to decay at infinity, which for the Dirichlet problem this implies that we have to change  the space $T^s$ for the subspace of functions such that $\langle u, 1 \rangle = 0$ (cf.~\cite{stephan1984augmented,JHP20}).
\end{remark}
\subsection{Elastic Wave Operators}
\label{sec:elastic}
In this case, one has that $\cP = \alpha \Delta + (\alpha + \beta) \nabla \nabla \cdot + \omega^2$. The  parameters $\alpha,\beta$ are called Lam\'e parameters\footnote{Typically, these are denoted $\mu, \lambda$ but we have changed this convention to $\alpha,\beta$ so as to avoid any confusion with solutions of Dirichlet and Neumann problems, $\bla$ and $\bmu$, respectively.}, with $\alpha >0$, $\alpha + \beta>0$, and $\omega>0$ is the angular frequency. We also define two associated wavenumbers $$
k_p^2 := \frac{\omega^2}{\alpha+2\beta}, \quad
k_s^2 := \frac{\omega^2}{\beta}.$$  
The co-normal trace corresponds to the traction operator, defined as: \begin{align*}
\mathcal{B}_{\bn} \mathbf{u} =
\alpha \mathbf{n} (\div(\mathbf{u}))
+ 2 \beta \partial_{\mathbf{n}} \mathbf{u}+
\beta  \mathbf{n}^\perp (\div(\mathbf{u}^\perp)),
\end{align*}
where for $\mathbf{v}=(v_1,v_2)$ we set $\mathbf{v}^\perp  = (v_2,-v_1)$. The standard  condition at infinity is called the {\em Kupradze radiation condition} \cite{Kupradze}.
We refer to \cite{kress1996} for uniqueness of the related BVPs. In this case, the fundamental solution is given by 
\begin{align}
\label{eq:efunsol}
\mathbf{G} (\mathbf{x}, \mathbf{y}) := 
\frac{i}{4 \beta} H^{(1)}_0(k_s d)\mathbf{I}+ \frac{i}{4\omega^2} \nabla_\mathbf{x} 
\nabla_\mathbf{x} \cdot \left( 
H^{(1)}_0(k_sd ) - H^{(1)}_0(k_p d)
\right),
\end{align}
where $d := \| \bx -\by\|$, and $\mathbf{I}$ denotes the identity matrix. Alternatively, following \cite{kress1996} this can be expressed as 
\begin{align*}
\mathbf{G}(\mathbf{x},\mathbf{y})  = G^1(d) \mathbf{I} + 
{G}^2(d)\mathbf{D}(\mathbf{x}-\mathbf{y}),
\end{align*}
where  $D(\mathbf{d}) = \dfrac{\mathbf{d} \mathbf{d}^t}{\|\mathbf{d}\|^2}$, and 
\begin{align*}
G^1(d) &:= \frac{i}{4\beta} H_{0}^{(1)}(k_s d) - \frac{i}{4\omega^2d}\left(k_s H_1^{(1)}(k_s d)- k_p H_1^{(1)}(k_p d)\right), \\
G^2(d) &:= \frac{i}{4\omega^2} \left( 
\frac{2k_s H^{(1)}_1(k_s d)-2k_p H^{(1)}_1(k_p d)}{d}+
k_p^2H^{(1)}_0(k_p d)- k_s^2H^{(1)}_0(k_s d)  
\right)
\end{align*}
Using the expansion of Hankel functions \cite[9.1.10 and 9.1.11]{abramowitz1965handbook}, we can express $G^1,G^2$ as 
$$
G^j(d) = R^j(d) +(\log d^2)J^j(d), \quad j=1,2.
$$
where $R^1,R^2$ are entire functions on the variable $d^2$, and 
\begin{align*}
J^1(d) &:= -\frac{J_0(k_s d)}{4\pi\beta}+ \frac{1}{4\pi\omega^2 d}\left(k_sJ_1(k_s d) - k_p J_1(k_p d)\right), \\
J^2(d) &:= -\frac{1}{4 \pi \omega^2} \left(
\frac{2k_s J_1(k_sd)-2k_p J_1(k_p d)}{d} +k_p^2J_0(k_p d) -k_s^2J_0(k_s d)
\right).
\end{align*}
Hence, from the series expansion of Bessel functions of zeroth and first order, we have that $J^1,J^2$ are entire functions in the $d^2$ variable, and also that
\begin{align*}
J^1(0) = -\frac{1}{4\pi\beta} -\frac{1}{8 \pi\omega^2}(k_p^2-k_s^2), \quad J^2(0) = 0.
\end{align*}
Thus, we can express the fundamental solution as \begin{align}\label{eq:elastdecomp}
	\mathbf{G}(\x,\y) 
	= 
	(\log d^2)\mathbf{J}
	+ 
	\mathbf{R},
\end{align}
where 
$\mathbf{J} = J^1(d) \mathbf{I} +J^2(d)\mathbf{D}(\x-\y)$,
and $\mathbf{R} = R^1(d) \mathbf{I} +R^2(d) \mathbf{D}(\x-\y)$.
The only difference with the canonical expression \eqref{eq:funsolgen}
is the presence of a factor $\mathbf{D}(\x-\y)$. 
Let us study the properties of this factor. 

\begin{lemma}
\label{lemma:Dmatrix}
Consider two arcs $\bm{r}$, and $\bp$ and define
the matrix function $\mathbf{D}_{\bm{r},\bp}$ as
\begin{align*}
	(\mathbf{D}_{\bm{r},\bp}(t,\tau))_{j,k} 
	\coloneqq
	\left\lbrace
	\begin{array}{cl}
		\dfrac{(r_j(t)-p_j(\tau))
		\cdot 
		(r_k(t)-p_k(\tau))}{d_{\bm{r},\bp}^2(t,\tau)},
		&
		\bm{r} \neq \bp\text{, or, } t \neq s \\
		\dfrac{r'_j(t) r'_k(\tau)}{\bm{r}'(t)\cdot \bm{r}'(\tau)},
		&
		\text{otherwise}
	\end{array}
	\right. ,
\end{align*} 
where $j,k \in {1,2}$, and also two compact sets $K^1,K^2 \subset \rgeoh{m}{\alpha}$ for some $m \in \IN$, and $\alpha \in [0,1]$. 
\begin{enumerate}
\item[(i)]
For $K^1 = K^2$, if we select $\delta$ as in Condition \ref{condition:Qself}, then it holds that
\begin{align*}
\bm{r} \in K^1 \mapsto ({D}_{\bm{r},\bm{r}})_{j,k} \in \cmspaceh{m-1}{\alpha}{(-1,1)\times(-1,1)}{\IC}, \quad j,k =1,2,
\end{align*} 
has a holomorphic extension in $K^1_\delta$.
\item[(ii)] 
For $K^1, K^2$ disjoint sets, if $\delta_1, \delta_2$ as in Condition  \ref{condition:dcross}, then the map
\begin{align*}
(\bm{r},\bp)  \in K^1 \times K^2 \mapsto ({D}_{\bm{r},\bp})_{j,k} \in \cmspaceh{m}{\alpha}{(-1,1)\times(-1,1)}{\IC}, \quad j,k =1,2,
\end{align*} 
has a holomorphic extension in $K^1_{\delta_1} \times K^2_{\delta_2}$.
\end{enumerate}
\end{lemma} 
\begin{proof}
For the first part, if $t \neq \tau$ we have that
\begin{align*}
(D_{\bm{r},\bm{r}}(t,\tau))_{j,k} = Q^{-1}_{\bm{r}}(t,\tau) \left( \frac{r_j(t)-r_j(\tau)}{t-\tau} \right) \cdot \left( \frac{r_k(t)-r_k(\tau)}{t-\tau} \right),
\end{align*}
where again $Q^{-1}_{\bm{r}}(t,\tau) = 1/Q_{\bm{r}}(t,\tau)$, and the results follow as in the proof of Lemma \ref{lemma:Qfun}.
The second part is direct from Lemma \ref{lemma:squaredistance} and elementary results of complex variable. 
\end{proof}

\
\begin{corollary}[Elastic case]
Consider the Dirichlet and Neumann BIEs (Problem \ref{prob:BIEs}) for the time-harmonic elastic kernel \eqref{eq:efunsol} for $\alpha >0$, $\alpha + \beta>0$, and $\omega>0$. Then, the arising BIOs and domain-to-solution maps are shape holomorphic.
\end{corollary}

\begin{proof}
We only need to ensure that the integral kernels can be expressed as in \ref{eq:G_disjoint_arcs}. The result for the Dirichlet problem follows directly by the decomposition of the fundamental solution \eqref{eq:elastdecomp}, Lemma \ref{lemma:Dmatrix} and also Remark \ref{remark:lessregularitykernel} for the $D_{\bm{r},\bm{r}}$ factor. 

For the Neumann problem, we use the following formula \cite[Equation 3.9]{tao2021},
\begin{align}
\label{eq:wstokeexp}
\begin{split}   
(W_{\Gamma_1,\hdots,\Gamma_M})_{i,j} \bu =  
\int_{-1}^1 \mathbf{G}_1(\bm{r}_i(t),\bm{r}_j(\tau))u_j(\tau) \text{d}s\\+
 \frac{\text{d}}{\text{d}t} \int_{-1}^1 \mathbf{G}_2(\bm{r}_i(t),\bm{r}_j(\tau)) \frac{\text{d} u_j(\tau)}{\text{d}s}\text{d}s+
\int_{-1}^{1} \mathbf{G}_3(\bm{r}_i(t),\bm{r}_j(\tau)) \frac{\text{d} u_j(\tau)}{\text{d}s}\text{d}s\\+ 
\frac{\text{d}}{\text{d}t} \int_{-1}^1 \mathbf{G}_2(\bm{r}_i(t),\bm{r}_j(\tau)) u_j(\tau)\text{d}s,
\end{split}
\end{align}
where the first kernel function is 
\begin{align*}
\begin{split}
&\mathbf{G}_1(\bm{r}_i(t),\bm{r}_j(\tau)) := \frac{i}{4}\left( \rho \omega^2(\bm{r}_i'^\perp(t)\bm{r}_j'^\perp(\tau)^T-\bm{r}_i'(t) \cdot \bm{r}_j'(\tau) \mathbf{I}) H^{(1)}_0(k_s d) \right. \\ 
&-\beta k_s^s(
\bm{r}_j'^\perp(\tau)\bm{r}_i'^\perp(t)^T 
-\bm{r}_i'^\perp(t)\bm{r}_j'^\perp(\tau)^T)H_0^{(1)}(k_sd) - \\ 
&\left. \rho\omega^2 (\bm{r}_i'^\perp(t)\bm{r}_j'^\perp(\tau)^T)H_0^{(1)}(k_p d)\right),
\end{split}
\end{align*}
which only has logarithmic singularities as it is composed of a combination zeroth-order Hankel functions. Hence, it corresponds to the term $\widetilde{G}$ in Maue's formula \eqref{eq:mauesrep}. The holomorphic extension of the corresponding BIO is direct since the kernel function is of the form described in Remark \ref{remark:lessregularitykernel}. 

The second kernel in \eqref{eq:wstokeexp} is 
\begin{align*}
\mathbf{G}_2(\bm{r}_i(t),\bm{r}_j(\tau))=
4\beta^2 \mathbf{A}\mathbf{G}(\bm{r}_i(t),\bm{r}_j(\tau)) \mathbf{A} + i \beta H_0^{(1)}(k_s d)\mathbf{I},
\end{align*}
where $\mathbf{A} = \begin{pmatrix}
0&-1\\1&0
\end{pmatrix}$. From the analysis of the weakly-singular BIO, we obtain the holomorphic extension of the corresponding integral operator. 

The third term in \eqref{eq:wstokeexp} is
\begin{align*}
\frac{i \beta \bm{r}_i'^\perp(t)(\bm{r}_i(t)-\bm{r}_j(\tau))}{2d}  \left(
k_s H^{(1)}_1(k_s d)  -
 k_p H^{(1)}_1(k_p d) 
\right)\mathbf{A},
\end{align*}
which can be shown not to have any singularities\footnote{This is almost the same as $G^1$ found in the analysis of the weakly-singular operator.}, but the structure still is the one described in Remark \ref{remark:lessregularitykernel}, only with $G_1 = 0$. The associated integral BIO in \eqref{eq:wstokeexp} can be seen as a map in $\mathcal{L}\left( \IU^s,\IW^s \right)$, where the range is in $\IW^s$ instead of $\IY^s$. This is due to the fact that if we start in $\IU^s$, by \eqref{eq:devprop} the derivative changes the argument in the operator to a function in $\IT^{s-1}$. Thus, by Corollary \ref{corollary:rfop}, we obtain the mentioned mapping property. Though operators with the mentioned mapping properties were not studied, by \eqref{eq:WsYs}, we can still consider that this operator lies in $\mathcal{L}\left(\IU^s,\IY^s \right)$, and hence it is compact in  $\mathcal{L}\left(\IU^s,\IY^{s-1} \right)$. The corresponding holomorphic extension then follows arguing as in the regular part of Theorem \ref{lemma:vcomp}. The final kernel function is \begin{align*}
\mathbf{G}_4(\bm{r}_i(t),\bm{r}_j(\tau)) = \frac{i\beta (\bm{r}_i(t)-\bm{r}^\perp(\tau))\bm{r}_j'^\perp(\tau)^T}{2d}
\left(
k_s H^{(1)}_1(k_s d)  -
k_p H^{(1)}_1(k_p d)
\right),
\end{align*}
which is also a regular kernel, whose structure is described in Remark \ref{remark:lessregularitykernel}, with $G_1 =0$. The mapping properties of the associated operator are not easily derived since, if $u_j \in \IU^s$, the evaluation of the BIO, discarding the derivative, would lie in $\IY^{s+1}$. Yet, we do not have a characterization of the derivative map in $\IY^{s+1}$. We circumvent this by using \eqref{eq:UsTs}, and so if $u_j \in \IU^s$ then $u_j \in \IT^s$, and hence the evaluation of the integral operator is in $\IW^{s+1}$. Thus, by \eqref{eq:devprop} the full operator is in $\mathcal{L}\left(\IU^s,\IY^s\right)$ and compact in $\mathcal{L}\left(\IU^s,\IY^{s-1}\right)$. Then, the holomoprphic extension follows as in the previous case. 
\end{proof}

\section{Conclusions and Future Work}
\label{sec:conclu}

We have shown a general framework for establishing parametric shape holomorphy of BIEs in two-dimensional space with multiple arcs. Though we have limited our findings to the case of homogeneous media, heterogenous coefficients and non-explicit fundamental solutions could also be addressed. Indeed, as long as the kernel can be decomposed as \eqref{eq:funsolgen} all results hold. Future work involves the application of these results in UQ and deep learning.


\bibliography{PHJ_22}
\bibliographystyle{siam}

\appendix

\section{Immersion of H{\"o}lder Spaces}
\label{appendix:lemma1}

{Herein, the symbol $\mathbf{n}$ will be used to denote a vector of two integer, not the normal of an arc. Also, we denote the non-normalized Fourier basis by $e_n(t) = \exp{(i n t) }$, for $n\in\mathbb{Z}$, and the bi-periodic basis as $e_{n,l}(t,\tau) = e_n(t)e_l(\tau)$, $n,l \in \mathbb{Z}$.

We begin by introducing the Sobolev-Slobodeckij norm for bi-periodic functions with domain  $[-\pi,\pi]\times[-\pi,\pi]$. Let $\mathbf{n} = (n_1,n_2) \in \IN_0^2$, and $\boldsymbol{\gamma} = (\gamma_1, \gamma_2) \in [0,1]^2$, we define 
\begin{align*}
\|g\|_{\boldsymbol{n},\boldsymbol{\gamma}}^2 
	= 
	\sum_{p \leq n_1} 
	\sum_{q \leq n_2} 
	\| \partial_t^{p} \partial_s^{q} g(t,\tau) \|^2_{L^2([-\pi,\pi]\times [-\pi,\pi])}
	+  
	|\partial_t^{n_1} \partial_{s}^{n_2} g(t,\tau)|_{\boldsymbol{\gamma}}^2,
\end{align*}
where the Sobolev-Slobodeckij semi-norm is defined as 
\begin{align*}
	|u|_{\boldsymbol{\gamma}}^2 
	= 
	\int_{-\pi}^{\pi}  
	\int_{-\pi}^{\pi} 
	\int_{-\pi}^{\pi}
	\int_{-\pi}^{\pi} 
	\frac{
		|u(x,y)-u(t,y)+u(t,\tau)-u(x,s)|^2
	}{
		\left(\sin\frac{|x-t|}{2}\right)^{1+2\gamma_1}\left(\sin\frac{|y-s|}{2}\right)^{1+2\gamma_2}
	}\text{d}x\text{d}y\text{d}t\text{d}s.
\end{align*}
 The semi-norm is generated by the following inner product:
\begin{align*}
\langle u , v \rangle{\boldsymbol{\gamma}} := \int_{-\pi}^{\pi}  \int_{-\pi}^{\pi} \int_{-\pi}^{\pi}\int_{-\pi}^{\pi} \frac{\widetilde\Delta u (x,t,y,s)\widetilde\Delta \overline{v} (x,t,y,s) }
{\left(\sin\frac{|x-t|}{2}\right)^{1+2\gamma_1}\left(\sin\frac{|y-s|}{2}\right)^{1+2\gamma_2}}\text{d}x\text{d}y\text{d}t\text{d}s, 
\end{align*}
where, the difference operator $\widetilde\Delta$ is given by
$$\widetilde\Delta u(x,t,y,s) := u(x,y)-u(t,y)+u(t,\tau)-u(x,s).$$
Following the uni-variate case \cite[Theorem 8.6]{kress2013linear}, we derive the following two results. 

\begin{lemma}
\label{lemma:fourierort}
For $(n_1,l_1), (n_2,l_2) \in \mathbb{Z}^2$ we have that: 
$$
\langle e_{n_1,l_1} , e_{n_2,l_2} \rangle_{\boldsymbol{\gamma}} = 16 \widetilde{\delta}_{n_1,n_2} \widetilde{\delta}_{l_1,l_2} S_{\gamma_1}(n_2) S_{\gamma_2}(l_2), 
$$
where $\widetilde{\delta}_{n,m} = 2\pi \delta_{n,m}$, for $n,m \in \mathbb{Z}$, and
$$
S_a(n) =\int_{0}^{\pi} \left(\sin\frac{n}{2} u\right)^2 
\left(\sin\frac{|u|}{2}\right)^{-1-2a} \text{d} u .
$$
\end{lemma}
\begin{proof}
We compute the inner product between two basis, $e_{n_1,l_1}$ and $e_{n_2,l_2}$.  By permuting integration variables, one has that 
$$\langle e_{n_1,l_1} , e_{n_2,l_2} \rangle{\boldsymbol{\gamma}}  =
4 \int_{-\pi}^{\pi}  \int_{-\pi}^{\pi} \int_{-\pi}^{\pi}\int_{-\pi}^{\pi} \frac{e_{n_1,l_1}(x,y)\widetilde\Delta e_{-n_2,-l_2} (x,t,y,s) }
{\left(\sin\frac{|x-t|}{2}\right)^{1+2\gamma_1}\left(\sin\frac{|y-s|}{2}\right)^{1+2\gamma_2}}\text{d}x\text{d}y\text{d}t\text{d}s.
$$
The right-hand side term is decomposed into the sum of four integrals defined as
\begin{align*}
I_1 &:= \int_{-\pi}^{\pi}\int_{-\pi}^{\pi} \int_{-\pi}^{\pi}\int_{-\pi}^{\pi} \frac{e_{n_1,l_1}(x,y)e_{-n_2,-l_2} (x,y) }
{\left(\sin\frac{|x-t|}{2}\right)^{1+2\gamma_1}\left(\sin\frac{|y-s|}{2}\right)^{1+2\gamma_2}}\text{d}x\text{d}y\text{d}t\text{d}s,\\
I_2 &:=  -\int_{-\pi}^{\pi}\int_{-\pi}^{\pi} \int_{-\pi}^{\pi}\int_{-\pi}^{\pi} \frac{e_{n_1,l_1}(x,y)e_{-n_2,-l_2} (t,y) }
{\left(\sin\frac{|x-t|}{2}\right)^{1+2\gamma_1}\left(\sin\frac{|y-s|}{2}\right)^{1+2\gamma_2}}\text{d}x\text{d}y\text{d}t\text{d}s,\\
I_3 &:= \int_{-\pi}^{\pi}\int_{-\pi}^{\pi} \int_{-\pi}^{\pi}\int_{-\pi}^{\pi} \frac{e_{n_1,l_1}(x,y)e_{-n_2,-l_2} (t,\tau) }
{\left(\sin\frac{|x-t|}{2}\right)^{1+2\gamma_1}\left(\sin\frac{|y-s|}{2}\right)^{1+2\gamma_2}}\text{d}x\text{d}y\text{d}t\text{d}s,\\
I_4 &:= -\int_{-\pi}^{\pi}\int_{-\pi}^{\pi} \int_{-\pi}^{\pi}\int_{-\pi}^{\pi} \frac{e_{n_1,l_1}(x,y)e_{-n_2,-l_2} (x,s) }
{\left(\sin\frac{|x-t|}{2}\right)^{1+2\gamma_1}\left(\sin\frac{|y-s|}{2}\right)^{1+2\gamma_2}}\text{d}x\text{d}y\text{d}t\text{d}s.
\end{align*}
We perform the change of variables $u = t-x$, and $v = s-y$ in the four integrals, and by the periodicity of the involving factors, we get 
\begin{align*}
I_1 &:= \widetilde{\delta}_{n_1,n_2}\widetilde{\delta}_{l_1,l_2}\int_{-\pi}^{\pi}
\left(\sin\frac{|u|}{2}\right)^{-1-2\gamma_1}\text{d}u
\int_{-\pi}^{\pi}  
\left(\sin\frac{|v|}{2}\right)^{-1-2\gamma_2}\text{d}v,\\
I_2 &:= -\widetilde{\delta}_{n_1,n_2}\widetilde{\delta}_{l_1,l_2}\int_{-\pi}^{\pi} e_{-n_2}(u)
\left(\sin\frac{|u|}{2}\right)^{-1-2\gamma_1}\text{d}u
\int_{-\pi}^{\pi}  
\left(\sin\frac{|v|}{2}\right)^{-1-2\gamma_2}\text{d}v,\\
I_3 &:= \widetilde{\delta}_{n_1,n_2}\widetilde{\delta}_{l_1,l_2}\int_{-\pi}^{\pi} e_{-n_2}(u)
\left(\sin\frac{|u|}{2}\right)^{-1-2\gamma_1}\text{d}u
\int_{-\pi}^{\pi}  e_{-l_2}(v)
\left(\sin\frac{|v|}{2}\right)^{-1-2\gamma_2}\text{d}v,\\
I_4 &:= -\widetilde{\delta}_{n_1,n_2}\widetilde{\delta}_{l_1,l_2}\int_{-\pi}^{\pi} 
\left(\sin\frac{|u|}{2}\right)^{-1-2\gamma_1}\text{d}u
\int_{-\pi}^{\pi}  e_{-l_2}(v)
\left(\sin\frac{|v|}{2}\right)^{-1-2\gamma_2}\text{d}v.
\end{align*}
Since the Fourier basis elements are $e_n(t) = \cos(nt) + i \sin (nt) $, onw can use the symmetries of  cosine and sine functions to obtain
$$
\int_{-\pi}^{\pi} (1-e_{-n_2}(u))
\left(\sin\frac{|u|}{2}\right)^{-1-2\gamma_1}\text{d}u = 
2 \int_{0}^{\pi} (1-\cos(n_2 u))
\left(\sin\frac{|u|}{2}\right)^{-1-2\gamma_1}\text{d}u.
$$
Furthermore, we invoke the double-angle formula for the cosine so as to arrive at
$$
\int_{-\pi}^{\pi} (1-e_{-n_2}(u))
\left(\sin\frac{|u|}{2}\right)^{-1-2\gamma_1}\text{d}u = 
4 \int_{0}^{\pi} \left(\sin\frac{n_2 u}{2}\right) ^2
\left(\sin\frac{|u|}{2}\right)^{-1-2\gamma_1}\text{d}u.
$$
Using this, we find
\begin{align*}
\begin{split}
&\langle e_{n_1,n_2}, e_{l_1,l_2} \rangle_{\boldsymbol{\gamma}}= 16 \widetilde{\delta}_{n_1,n_2} \widetilde{\delta}_{l_1,l_2}\\ &\int_{0}^{\pi} \left(\sin\frac{n_2 u}{2}\right) ^2
\left(\sin\frac{|u|}{2}\right)^{-1-2\gamma_1}\text{d}u 
\int_{0}^{\pi} \left(\sin\frac{l_2 v}{2}\right) ^2
\left(\sin\frac{|v|}{2}\right)^{-1-2\gamma_2}\text{d}u.
\end{split}
\end{align*}
which is equivalent to the statement of the lemma.
\end{proof}

\begin{lemma}
\label{lemma:purefractionalbound}
Let $\boldsymbol{\gamma} =(\gamma_1,\gamma_2) \in [0,1]^2$, and $\varrho$ a bi-periodic function in $[-\pi,\pi]\times [-\pi,\pi]$,  then $$
\| \varrho\|_{\gamma_1,\gamma_2}^2 \lesssim \| \varrho\|^2_{L^2([-\pi,\pi]\times [-\pi,\pi])} + |\varrho|^2_{\boldsymbol{\gamma}}
,$$
i.e.~the Sobolev norm for two pure fractional orders $\gamma_1, \gamma_2$ of a bi-periodic function is bounded by the Sobolev-Slobodeckij norm of order $(0,0),(\gamma_1,\gamma_2)$.
\end{lemma}

\begin{proof}
We consider a function $\varrho$ expanded in terms of bi-periodic Fourier basis functions. By Lemma \ref{lemma:fourierort}, it holds that 
\begin{align}
|\varrho |_{\boldsymbol{\gamma}}^2 = \langle \varrho ,\varrho \rangle_{\boldsymbol{\gamma}} \cong \sum_{n=-\infty}^\infty \sum_{l=-\infty}^{\infty} |\widetilde{\varrho}_{n,l}|^2 S_{\gamma_1}(n) S_{\gamma_2}(l).
\end{align}
From \cite[(8.8)]{kress2013linear}, we have that  the function $S_a(n)$ from Lemma \ref{lemma:fourierort}, behave as $S_a(n) \cong (n^2)^a$, for $a \in [0,1]$. Consequently, 
\begin{align}
|\varrho |_{\boldsymbol{\gamma}}^2 \cong \sum_{n=-\infty}^\infty \sum_{l=-\infty}^{\infty} |\widetilde{\varrho}_{n,l}|^2 (n^2)^{\gamma_1} (l^2)^{\gamma_2}.
\end{align}
We can now use the well-known inequality  $(1+n^2)^{\gamma_1}  \leq (n^2)^{\gamma_1} +1$---analogously for $l$ and $\gamma_2$--, so that
 \begin{align*}
 \|\varrho\|_{\gamma_1,\gamma_2}^2 =
\sum_{n= -\infty}^\infty \sum_{l = -\infty}^\infty (1+n^2)^{\gamma_1} (1+l^2)^{\gamma_2}|\widetilde{\varrho}_{n,l}|^2 \lesssim
 |\varrho|^2_{\boldsymbol{\gamma}} + \| \varrho\|^2_{L^2([-\pi,\pi]\times [-\pi,\pi])}   \end{align*}
as stated.
\end{proof}

\subsection{Proof of Lemma \ref{lemma:cmdecay2}}

We will first show that  $\|g\|_{s_1,s_2}^2$ can be bounded by the Sobolev-Slobodeckij norm $ \|g\|_{\mathbf{n},\boldsymbol{\gamma}}^2$, with $\bf{n} = ([s_1],[s_2])$, and $\boldsymbol{\gamma} = (\{s_1\},\{s_2\})$, where $[s_1], [s_2]$ denote the integer parts of $s_1, s_2$ respectively and $\{s_1\}, \{s_2\}$ are the corresponding fractional parts. By definition, one has that 
$$
\|g\|_{s_1,s_2}^2  =  \sum_{n = -\infty}^\infty \sum_{l = -\infty}^\infty (1+n^2)^{\{s_1\}} (1+l^2)^{\{s_2\}} (1+n^2)^{[s_1]} (1+l^2)^{[s_2]} \left\vert\widetilde{g}_{n,l}\right\vert^2
$$
Using the inequality $(1+n^2)^s \lesssim (n^2)^s+1$, we get
\begin{align*}
\begin{split}
    |g\|_{s_1,s_2}^2 & \lesssim \sum_{n = -\infty}^\infty \sum_{l = -\infty}^\infty (1+n^2)^{\{s_1\}} (1+l^2)^{\{s_2\}} (1+n^{2[s_1]}) (1+l^{2[s_2]}) \left\vert\widetilde{g}_{n,l}\right\vert^2\\
 &\lesssim \sum_{p=0}^{[s_1]} \sum_{q=0}^{[s_2]}\sum_{n = -\infty}^\infty \sum_{l = -\infty}^\infty (1+n^2)^{\{s_1\}} (1+l^2)^{\{s_2\}} (1+n^{2p}) (1+l^{2q}) \left\vert\widetilde{g}_{n,l}\right\vert^2
\end{split}
\\
&\lesssim
\sum_{p=0}^{[s_1]} \sum_{q=0}^{[s_2]}\sum_{n = -\infty}^\infty \sum_{l = -\infty}^\infty (1+n^2)^{\{s_1\}} (1+l^2)^{\{s_2\}} (1+n^{2p}l^{2q})  \left\vert\widetilde{g}_{n,l}\right\vert^2
\end{align*}
We notice that the Fourier coefficients of $\partial_t^{p} \partial_s^{q} g$ are in fact $(i n)^{p} (i l)^{q} \widetilde{g}_{n,l}$. Hence, we obtain that
\begin{align}
\label{eq:intdevs}
\begin{split}
\|g\|_{s_1,s_2}^2 
\lesssim
\sum_{p=0}^{[s_1]} \sum_{q=0}^{[s_2]}\sum_{n = -\infty}^\infty \sum_{l = -\infty}^\infty (1+n^2)^{\{s_1\}} (1+l^2)^{\{s_2\}}  \left\vert\left(\widetilde{\partial_t^{[p]} \partial_s^{[q]}g(t,\tau)}\right)_{n,l}\right\vert^2+ \\
 \sum_{n = -\infty}^\infty \sum_{l = -\infty}^\infty (1+n^2)^{\{s_1\}} (1+l^2)^{\{s_2\}} \left\vert\widetilde{g}_{n,l}\right\vert^2
,
\end{split}
\end{align}
 From the last inequality, we notice that for  $\{s_1\}= \{s_2\} =0$ so that
$$
\|g\|_{s_1,s_2}^2 \lesssim \sum_{p=0}^{[s_1]} \sum_{q=0}^{[s_2]}\sum_{n = -\infty}^\infty \sum_{l = -\infty}^\infty  \left\vert\left(\widetilde{\partial_t^{[p]} \partial_s^{[q]}g(t,\tau)}\right)_{n,l}\right\vert^2+
 \sum_{n = -\infty}^\infty \sum_{l = -\infty}^\infty  \left\vert\widetilde{g}_{n,l}\right\vert^2.
$$
Then, by Parseval's identity we obtain 
$$
\|g\|^2_{s_1,s_2} \lesssim \|g\|^2_{([s_1],[s_2]), (0,0)}
$$
In any other case, $\{s_1\}>0$ or $\{s_2\}>0$, for every $p \in \{0, \hdots, [s_1]\}$ and $q \in \{0, \hdots,[s_2]\}$ we define $\varrho^{p,q} = \partial_t^{p} \partial_s^{q}g(t,\tau)$, and from \eqref{eq:intdevs} we see that we only need to show that 
\begin{align*}
 \sum_{n = -\infty}^\infty \sum_{l = -\infty}^\infty (1+n)^{\{s_1\}} (1+l^2)^{\{s_2\}} \left\vert\widetilde{\varrho^{p,q}}_{n,l}\right\vert^2 \lesssim \| \varrho \|^2_{(0,0),(\{s_1\},\{s_2\})},
\end{align*}
which holds by Lemma \ref{lemma:purefractionalbound}  with $\varrho := \varrho^{[s_1],[s_2]}$. We conclude that 
\begin{align}
\label{eq:gbound1}
    \|g\|_{s_1,s_2}^2 \lesssim
    \|g\|_{([s_1],[s_2]),(\{s_1\},\{s_2\})}^2.
\end{align}
To finish the proof we will bound the general norm $\|g\|_{\mathbf{n},\boldsymbol{\gamma}}^2$, in terms of a $\cmspaceh{m}{\alpha}{[-\pi,\pi]\times[-\pi,\pi]}{\IC}$-norm with appropiate $m \in \IN_0$, and $\alpha \in [0,1]$.
First, from the definition of  H\"older norms, it holds that
\begin{align*}
 \|g\|^2_{\mathbf{n},(0,0)} \lesssim \|g\|_{\cmspaceh{m}{\alpha}{[-\pi,\pi]\times[-\pi,\pi]}{\IC}},
\end{align*}
for $m+\alpha \geq n_1 + n_2$. For the purely fractional case, we have that for any $a, b >0$ such that $a + b = 1$
\begin{align*}
\begin{split}
|g|_{(\gamma_1, \gamma_2)}^2 &\lesssim \|g\|_{\cmspaceh{0}{\alpha}{[-\pi,\pi]\times [-\pi,\pi]}{\IC}}^2\\&\cj{\times}\left( \int_{-\pi}^{\pi}\int_{-\pi}^{\pi} \frac{|x-t|^{2 a \alpha}}{\left(\sin\frac{|x-t|}{2} \right)^{1+2\gamma_1}} \text{d}t \text{d}x\right)
\left( \int_{-\pi}^{\pi}\int_{-\pi}^{\pi} \frac{|y-s|^{2 b \alpha }}{\left(\sin\frac{|y-s|}{2} \right)^{1+2\gamma_2}} \text{d}s \text{d}y\right),
\end{split}
\end{align*}
and thus, the right-hand side is finite only if $\alpha > \gamma_1 + \gamma_2$. However, this condition cannot be used whenever $\gamma_1 +\gamma_2 \geq 1$. For the latter, we assume that $m \geq 1$, and by the mean value theorem, we see that 
\begin{align*}
\begin{split}
|g(x,y)-g(t,y)+g(t,\tau)-g(x,s)|
&=  \left\vert \int_{t}^x \partial_1 g(\lambda,y) - \partial_1 g(\lambda,s)\text{d} \lambda \right\vert \\ &=|x-t|
 \left\vert\partial_1 g(\xi,y)-\partial_t g(\xi,s)\right\vert\\ &\leq |x-t||y-s|^\alpha \|g \|_{\cmspaceh{1}{\alpha}{[-\pi,\pi]\times[-\pi,\pi]}{\IC}}.
\end{split}
\end{align*}
Similarly, by reordering terms we conclude that
$$
|g(x,y)-g(t,y)+g(t,\tau)-g(x,s)|\leq |y-s||x-t|^\alpha \|g \|_{\cmspaceh{1}{\alpha}{[-\pi,\pi]\times[-\pi,\pi]}{\IC}}.
$$
Thus, for every $a \geq 0, b \geq 0$ such that $a + b = 1$ we have 
\begin{align*}
\begin{split}
|g|_{\boldsymbol{\gamma}}^2 &\lesssim \|g \|^2_{\cmspaceh{1}{\alpha}{[-\pi,\pi]\times[-\pi,\pi]}{\IC}} \\
&\times\int_{-\pi}^{\pi}  \int_{-\pi}^{\pi} \int_{-\pi}^{\pi}\int_{-\pi}^{\pi} \frac{ (|x-t|^{a + \alpha b}  
|y-s|^{b + \alpha a}
)^2}
{\left(\sin\frac{|x-t|}{2}\right)^{1+2\gamma_1}\left(\sin\frac{|y-s|}{2}\right)^{1+2\gamma_2}}\text{d}x \text{d}y \text{d}t \text{d}s
\end{split},
\end{align*}
and one can easily see that integrals in the right hand side are finite if $\alpha > \gamma_1 + \gamma_2 -1$. 
Finally, we conclude that 
\begin{align*}
\|g\|_{\mathbf{n},\boldsymbol{\gamma}}^2 &\lesssim \|g\|^2_{\cmspaceh{m}{\alpha}{[-\pi,\pi]\times[-\pi,\pi]}{\IC}}\\&\times\begin{cases}
 \text{  for $m_1+m_2+\alpha_1 + \alpha_2 < n+\gamma$ if $\alpha_1 + \alpha_2 < 1$},\\
 \text{ for $(n_1+n_2+1)+(\gamma_1+\gamma_2 -1)< m+\alpha$, if $\gamma_1 + \gamma_2 \geq 1$}.
\end{cases}
\end{align*}
which is equivalent to 
$$
\|g\|_{\mathbf{n},\boldsymbol{\gamma}}^2 \lesssim \|g\|^2_{\cmspaceh{m}{\alpha}{[-\pi,\pi]\times[-\pi,\pi]}{\IC}},
$$
if $n_1+n_2+\gamma_1+\gamma_2 < m+\alpha$. 
The statement of the lemma then follows from this last bound and \eqref{eq:gbound1}.
\qed

\section{Proofs Lemmas \ref{lemma:squaredistance} and \ref{lemma:Qfun}}
\label{appendix:funext}
We present proofs for Lemmas \ref{lemma:squaredistance} and \ref{lemma:Qfun}. For the latter, we will need the following auxiliary result.

\begin{lemma}
\label{lemma:dwelldef}
Let $m \in \IN$, $\alpha \in [0,1]$ and $K \subset \rgeoh{m}{\alpha}$ be a compact set of $\cmspaceh{m}{\alpha}{(-1,1)}{\IR^2}$.
\begin{enumerate}
	\item[(i)]
	It holds that
	\begin{align}
		\inf_{\bm{r} \in K } \inf_{t \in (-1,1)} \| \bm{r}'(t) \|>0  
		\quad 
		\text{and} 
		\quad  
		\sup_{\bm{r} \in K} \sup_{t \in (-1,1)} \| \bm{r}'(t)\| < \infty.
	\end{align}
	\item[(ii)]
	There exists $\delta >0 $ such that  
	\begin{equation}
		\inf_{\bm{r} \in K_\delta} \inf_{t \in (-1,1)}
		\mathfrak{Re}(\{
			(\bm{r}'(t) \cdot \bm{r}'(t))
		\}
		>0.
	\end{equation}
		i.e.~
	there exists $\delta>0$ fulfilling Condition \ref{condition:Qself}.
\end{enumerate}
\end{lemma}
\begin{proof}
The first part follows from the continuity of
\begin{equation}
	\mathcal{I}(\bm{r})
	\coloneqq
	\inf_{t \in (-1,1)} 
	\| \bm{r}'(t)\| 
	\quad
	\text{and} 
	\quad
	\mathcal{S}(\bm{r}) 
	\coloneqq 
	\sup_{t \in (-1,1)} \| \bm{r}'(t)\|
\end{equation}
in $\rgeoh{m}{\alpha}$.
Indeed, if for each $\bm{r} \in \rgeoh{m}{\alpha}$
we consider any $\bp \in \cmspaceh{m}{\alpha}{(-1,1)}{\IR^2}$
such that  $\| \bm{r} -\bp \|_{\cmspaceh{m}{\alpha}{(-1,1)}{\IR^2}}< \epsilon$, one has that
\begin{align}
	\snorm{
	\mathcal{I}(\bm{r})
	-
	\mathcal{I}(\bp)}
	&
	\leq
	\inf_{t \in (-1,1)} 
	\snorm{
	\norm{\bm{r}'(t)}
	-
	\norm{\bp'(t)}
	} \\
	&
	\leq
	\inf_{t \in (-1,1)} 
	\norm{\bm{r}'(t)-\bp'(t)}
	\leq
	\| \bm{r} -\bp \|_{\cmspaceh{m}{\alpha}{(-1,1)}{\IR^2}}< \epsilon.
\end{align}
Thus, one concludes that the map 
$\bm{r} \mapsto \mathcal{I}(\bm{r})$ is continuous,
and then the infimum in $K$ is achieved since $K$ is compact.
The supremum case follows similarly. 

For the second part we set
\begin{equation}
	\mathcal{I}_1
	=
	\inf_{\bm{r} \in K } 
	\mathcal{I}(\bm{r})
	\quad
	\text{and}
	\quad
	\mathcal{S}_1 
	=
	\sup_{\bm{r} \in K } \mathcal{S}(\bm{r}),
\end{equation}
and consider any element $\bm{r} \in K_\delta$. Then, there is $\bp \in K$ such that $\| \bm{r} -\bp \|_{\cmspaceh{m}{\alpha}{(-1,1)}{\IR^2}} < \delta$, and it holds that
\begin{align*}
\bm{r}' \cdot \bm{r}' = \|\bm{r}'\|^2+ 2(\bm{r}' -\bp')\cdot \bp' +(\bm{r}'-\bp')\cdot(\bm{r}'-\bp') 
\end{align*}
Therefore,
\begin{align*}
\mathfrak{Re}(\bm{r}' \cdot \bm{r}') \geq \mathcal{I}_1^2 - 2\mathcal{S}_1\delta -\delta^2, \end{align*}
and the result then follows by selecting $\delta < \sqrt{\mathcal{I}_1^2+\mathcal{S}_1^2}-\mathcal{S}_1$.
\end{proof} 

\subsection{Proof of Lemma \ref{lemma:squaredistance}}
\label{proof:lem1}
First, we prove that the map 
$$\bm{r},\bp \in \cmspaceh{m}{\alpha}{(-1,1)}{\IC^2}  \times \cmspaceh{m}{\alpha}{(-1,1)}{\IC^2}  \rightarrow d_{\bm{r},\bp}^2 \in \cmspaceh{m}{\alpha}{(-1,1)}{\IC}$$ has a Fr\'echet derivative at every point. Consequently, the holomorphic extension in compact sets of  $\rgeoh{m}{\alpha} \times \rgeoh{m}{\alpha}$ follows directly. We will limit ourself to the cases where $m\geq 1$, since the result is formulated in terms of $\rgeoh{m}{\alpha}$ and by convention $\rgeoh{0}{\alpha} =\emptyset$.

Consider two arbitrary functions $\bm{r},\bp \in \cmspaceh{m}{\alpha}{(-1,1)}{\IC^2}$. By definition, we have that $$d_{\bm{r},\bp}^2(t,\tau) \in \cmspaceh{m}{\alpha}{\iinterv}{\IC}.$$
For a pair $\bv^1,\bv^2 \in \cmspaceh{m}{\alpha}{(-1,1)}{\IC^2}$, we define
$$Dd_{\bm{r},\bp}^2[\bv^1,\bv^2](t,\tau) := 2(\bm{r}(t)-\bp(\tau))\cdot \left(\bv^1(t)-\bv^2(\tau)\right),$$
which is linear in $\bv^1,\bv^2$ and lies in $\cmspaceh{m}{\alpha}{\iinterv}{\IC}$. Finally, we notice that 
$$
d^2_{\bm{r}+\bv^1,\bp+\bv^2}(t,\tau) -d^2_{\bm{r},\bp}(t,\tau) -Dd_{\bm{r},\bp}^2[\bv^1,\bv^2](t,\tau) = d^2_{\bv^1,\bv^2}(t,\tau).
$$
By the product derivation rule, we conclude that 
\begin{align*}
\|d^2_{\bv^1,\bv^2}(t,\tau)\|_{\cmspaceh{m}{\alpha}{\iinterv}{\IC}} \lesssim& 2\|\bv^1\|_{\cmspaceh{m}{\alpha}{(-1,1)}{\IC^2}} \|\bv^2\|_{\cmspaceh{m}{\alpha}{(-1,1)}{\IC^2}}  \\ &+\|\bv^1\|_{\cmspaceh{m}{\alpha}{(-1,1)}{\IC^2}}^2+\|\bv^2\|_{\cmspaceh{m}{\alpha}{(-1,1)}{\IC^2}}^2\\
\lesssim &
\|\bv^1\|_{\cmspaceh{m}{\alpha}{(-1,1)}{\IC^2}}^2+\|\bv^2\|_{\cmspaceh{m}{\alpha}{(-1,1)}{\IC^2}}^2.
\end{align*} 
Hence, $Dd_{\bm{r},\bp}^2[\bv^1,\bv^2]$ is in fact the Fr\'echet derivative of $d_{\bm{r},\bp}^2$. Since $\bm{r},\bp$ are arbitrary, the distance function is holomorphic in an arbitrary open set of 
$\cmspaceh{m}{\alpha}{(-1,1)}{\IC^2} \times \cmspaceh{m}{\alpha}{(-1,1)}{\IC^2} $.

For the final part, notice that since the sets are admissible, it holds that 
\begin{align*}
I_d = \inf_{(\bm{r},\bp) \in K^1 \times K^2} \inf_{(t,\tau) \in (-1,1)\times(-1,1)}
 \| \bm{r}(t) - \bp(\tau) \| >0.
\end{align*}
Using this along with the fact that for every $\delta_1, \delta_2$, $K^1_{\delta_1}$, and $K^2_{\delta_2}$ are bounded, we conclude that there exist $\delta_1, \delta_2$ satisfying Condition \ref{condition:dcross}. Furthermore, for every $\bm{r},\bp  \in K^1_{\delta_1} \times K^2_{\delta_2}$, we have that
$$
\mathfrak{Re}(d_{\bm{r},\bp}^2(t,\tau)) \geq \mathcal{I}_d^2 -2 \mathcal{S}_d(\delta_1 + \delta_2) - ( \delta_1 + \delta_2)^2,$$
where $\mathcal{S}_d$, is defined as in Condition \ref{condition:dcross}. Lastly, by considering the restriction on $\delta_1, \delta_2$ due to Condition \ref{condition:dcross}, we directly have  $$\mathfrak{Re}(d_{\bm{r},\bp}^2(t,\tau))>0.$$
\subsection{Proof of Lemma \ref{lemma:Qfun}}
\label{proof:lem2}
Consider $\bm{r} \in \cmspaceh{m}{\alpha}{(-1,1)}{\IC^2}$, with $m \geq 1$, for the same reasons required for the proof of Lemma \ref{proof:lem1}. By Taylor expansion, it is immediate that  $$Q_{\bm{r}} \in \cmspaceh{m-1}{\alpha}{(-1,1)\times(-1,1)}{\IC}.$$  
For $\bv \in \cmspaceh{m}{\alpha}{(-1,1)}{\IC^2}$, let us define
\begin{align*}
DQ_{\bm{r}}[\bv](t,\tau) = \frac{2 (\bm{r}(t)-\bm{r}(\tau))\cdot \left(\bv(t)-\bv(\tau)\right)}{(t-s)^2}, 
\end{align*}
with the continuous extension for $t=s$, given by the Taylor expansions of $\bm{r}$ and $\bv$. It is clear that $DQ_{\bm{r}}[\bv]$ is linear in the $\bv$ variable and also $DQ_{\bm{r}}[\bv] \in \cmspaceh{m-1}{\alpha}{(-1,1)\times(-1,1)}{\IC}$. We also have that 
\begin{align*}
d_{\bm{r} +\bv}^2(t,\tau) = d^2_{\bm{r}}(t,\tau) + 2 (\bm{r}(t) -\bm{r}(\tau))\cdot (\bv(t)-\bv(\tau)) + d_{\bv}^2(t,\tau).
\end{align*}
Therefore, one has that
\begin{align*}
Q_{\bm{r} +\bv}(t,\tau) = Q_{\bm{r}}(t,\tau) + DQ_{\bm{r}}[v](t,\tau) + Q_{\bv}(t,\tau).
\end{align*}
Arguing as in the proof of Lemma \ref{lemma:squaredistance}, 
we have that 
$$\| Q_{\bv}\|_{\cmspaceh{m-1}{\alpha}{\iinterv}{\IC^2}} \lesssim \| \bv\|_{\cmspaceh{m}{\alpha}{(-1,1)}{\IC^2}}^2,$$ and one concludes that there is a Fr\'echet derivative everywhere. Thus, the function has a holomorphic extension in $\rgeoh{m}{\alpha}$.

Now, we show that the real part is strictly positive. Using the Taylor expansion once again, we have that
\begin{align*}
d_{\bm{r}}^2(t,\tau) = (t-s)^2 \left(\int_{0}^1 \bm{r}'(t+\delta(s-t))\text{d}\delta \right) \cdot \left(\int_{0}^1 \bm{r}'(t+\delta(s-t))\text{d}\delta \right).
\end{align*}
Consequently, we can write
$$
Q_{\bm{r}}(t,\tau) = \left(\int_{0}^1 \bm{r}'(t+\delta(s-t))\text{d}\delta \right) \cdot \left(\int_{0}^1 \bm{r}'(t+\delta(s-t))\text{d}\delta \right).
$$
The result then follows from the mean value theorem and selecting $\delta$ as in Lemma \ref{lemma:dwelldef}. The results for $Q^{-1}_{\bm{r}}$ follows using the fact that the function $z^{-1}$ is holomorphic away from zero.

\end{document}